\documentclass[reqno]{amsart}

\usepackage{amsmath}
\usepackage{amsthm}
\usepackage{amsfonts}
\usepackage{amssymb}
\usepackage{latexsym}
\usepackage{amscd}
\usepackage{stmaryrd}
\usepackage{color}
\usepackage[all]{xypic}
\usepackage{epsfig}
\usepackage{graphics}
\usepackage{ifthen}
\usepackage{varioref}

\theoremstyle{plain}
\newtheorem{thm}{Theorem}[section]
\newtheorem{cor}[thm]{Corollary}
\newtheorem{lem}[thm]{Lemma}
\newtheorem{prop}[thm]{Proposition}
\newtheorem{rem}[thm]{Remark}
\newtheorem{dfn}[thm]{Definition}
\newtheorem{conj}[thm]{Conjecture}

\numberwithin{equation}{section}

\definecolor{darkgreen}{rgb}{0.0625,0.64,0.0625}
\usepackage[
            pdfstartview=FitH,
            CJKbookmarks=true,
            bookmarksnumbered=true,
            bookmarksopen=true,
            colorlinks,
            pdfborder=001,
            linkcolor=blue,
            anchorcolor=green,
            citecolor=red
            ]{hyperref}

\newfont{\scyr}{wncyr10 scaled 550}
\def\shuffle{\,\mbox{\bf \scyr X}\,}

\def\sshuffle{\overset{\scalebox{1.4}{\_}\hspace{0.05cm}}{\shuffle}}
\def\sast{\overset{\scalebox{1.4}{\_}\hspace{0.05cm}}{\ast}}

\def\proof{\noindent {\bf Proof.\;}}

\def\wt{\operatorname{wt}}
\def\dep{\operatorname{dep}}
\def\height{\operatorname{ht}}
\def\Li{\operatorname{Li}}
\def\reg{\operatorname{reg}}
\def\Gal{\operatorname{Gal}}

\allowdisplaybreaks

\begin{document}

\title{Some relations deduced from regularized double shuffle relations of multiple zeta values}

\date{\today\thanks{The first author is supported by the National Natural Science Foundation of
China (Grant No. 11471245).} }

\author{Zhonghua Li \quad and \quad Chen Qin}

\address{School of Mathematical Sciences, Tongji University, No. 1239 Siping Road,
Shanghai 200092, China}

\email{zhonghua\_li@tongji.edu.cn}

\address{School of Mathematical Sciences, Tongji University, No. 1239 Siping Road,
Shanghai 200092, China}

\email{2014chen\_qin@tongji.edu.cn}



\begin{abstract}
It is conjectured that the regularized double shuffle relations give all algebraic relations among the multiple zeta values, and hence all other algebraic relations should be deduced from the regularized double shuffle relations. In this paper, we provide as many as the relations which can be derived from the regularized double shuffle relations, for example, the weighted sum formula of L. Guo and B. Xie, some evaluation formulas with even arguments and the restricted sum formulas of M. E. Hoffman and their generalizations.
\end{abstract}

\maketitle


\tableofcontents


\section{Introduction}\label{Sec:Intro}

For a finite sequence $\mathbf{k}=(k_1,\ldots,k_n)$ of positive integers, we define its weight, depth and height by
$$\wt(\mathbf{k})=k_1+\cdots+k_n,\quad \dep(\mathbf{k})=n,\quad \height(\mathbf{k})=\sharp\{l\mid 1\leqslant l\leqslant n,k_l\geqslant 2\},$$
respectively. If $k_1>1$, we call $\mathbf{k}$ an admissible index. For such an admissible index $\mathbf{k}=(k_1,\ldots,k_n)$, the multiple zeta value $\zeta(\mathbf{k})$ and the multiple zeta-star value $\zeta^{\star}(\mathbf{k})$ are real numbers defined respectively by the following convergent series
\begin{align*}
&\zeta(\mathbf{k})=\zeta(k_1,\ldots,k_n)=\sum\limits_{m_1>\cdots>m_n\geqslant 1}\frac{1}{m_1^{k_1}\cdots m_n^{k_n}},\\
&\zeta^{\star}(\mathbf{k})=\zeta^{\star}(k_1,\ldots,k_n)=\sum\limits_{m_1\geqslant\cdots\geqslant m_n\geqslant 1}\frac{1}{m_1^{k_1}\cdots m_n^{k_n}}.
\end{align*}
When $n=1$, both cases degenerate to the Riemann zeta values, which are the special values of the Riemann zeta function at the positive integer arguments. The study of these numbers can be traced back to L. Euler \cite{Euler}, while the systematic research began with M. E. Hoffman's paper \cite{Hoffman1992} and D. Zagier's paper \cite{Zagier}.

It is easy to see that every multiple zeta-star value can be written as a $\mathbb{Z}$-linear combination of multiple zeta values and vice versa as shown in the following simple example
\begin{align*}
\zeta^{\star}(2,1)=\sum\limits_{m\geqslant n\geqslant 1}\frac{1}{m^2n}=\sum\limits_{m>n\geqslant 1}\frac{1}{m^2n}+\sum\limits_{m=n\geqslant 1}\frac{1}{m^2n}=\zeta(2,1)+\zeta(3).
\end{align*}
Hence the $\mathbb{Q}$-vector space $\mathcal{Z}$ spanned by all multiple zeta values coincides with that spanned by all multiple zeta-star values. One important thing in the study of multiple zeta values is to understand the algebraic structure of the space $\mathcal{Z}$. A recent theorem of F. Brown \cite{Brown} states that all periods of mixed Tate motives unramified over $\mathbb{Z}$ are $\mathbb{Q}\left[\frac{1}{2\pi i}\right]$-linear combinations of multiple zeta values, and the multiple zeta values indexed by $2$ and $3$ are linear generators of the $\mathbb{Q}$-vector space $\mathcal{Z}$. As a consequence, one obtains an upper bound for the dimension of the subspace $\mathcal{Z}_k$, which is spanned by all multiple zeta values of weight $k$. More specifically, we have
$$\dim_{\mathbb{Q}} \mathcal{Z}_k\leqslant d_k,$$
where the sequence $\{d_k\}$ is defined by $d_0=1$, $d_1=0$, $d_2=1$ and $d_k=d_{k-2}+d_{k-3}$ for all integers $k\geqslant 3$. This result was first proved by A. B. Goncharov \cite{Goncharov} and independently  by T. Terasoma \cite{Terasoma}. The well-known dimension conjecture claims that the dimension of the vector space $\mathcal{Z}_k$ is exactly $d_k$, which was proposed by D. Zagier \cite{Zagier}.

Since $d_k$ is far less than $2^{k-2}$, the number of the admissible indices of weight $k$, there must be a lot of linear (algebraic) relations among multiple zeta values. Here we recall the regularized double shuffle relations from \cite{Ihara-Kaneko-Zagier,Racinet}, which are conjectured to give all algebraic relations of multiple zeta values.

As a subset of the field $\mathbb{R}$ of real numbers, $\mathcal{Z}$ is not only a subspace but also a subalgebra, which is because a product of two multiple zeta values is a sum of multiple zeta values. There are two different ways to multiple two multiple zeta values. The first one uses the infinite series representations as displaying in the following example:
\begin{align*}
&\zeta(2)\zeta(2)=\left(\sum\limits_{m=1}^\infty\frac{1}{m^2}\right)\left(\sum\limits_{n=1}^\infty\frac{1}{n^2}\right)\\
=&\left(\sum\limits_{m>n\geqslant 1}+\sum\limits_{n>m\geqslant 1}+\sum\limits_{m=n\geqslant 1}\right)\frac{1}{m^2n^2}\\
=&2\zeta(2,2)+\zeta(4).
\end{align*}
Such type products are called the stuffle products. The second one uses the so called Drinfel'd iterated integral representations. In fact, for an admissible index $\mathbf{k}=(k_1,\ldots,k_n)$, we have
$$\zeta(\mathbf{k})=\int\limits_{1>t_1>\cdots>t_k>0}f_1(t_1)\cdots f_k(t_k)dt_1\cdots dt_k$$
with
\begin{align}
f_i(t)=\begin{cases}
\frac{1}{1-t} & \text{if\;} i=k_1,k_1+k_2,\ldots,k_1+\cdots+k_n=k,\\
\frac{1}{t} & \text{otherwise}.
\end{cases}
\label{Eq:One-Forms}
\end{align}
Using the iterated integral representations, we have, for example,
\begin{align*}
&\zeta(2)\zeta(2)=\left(\;\int\limits_{1>t_1>t_2>0}\frac{dt_1dt_2}{t_1(1-t_2)}\right)\left(\;\int\limits_{1>s_1>s_2>0}\frac{ds_1ds_2}{s_1(1-s_2)}\right)\\
=&\left(\;\int\limits_{1>t_1>t_2>s_1>s_2>0}+\int\limits_{1>t_1>s_1>t_2>s_2>0}+\int\limits_{1>t_1>s_1>s_2>t_2>0}+\int\limits_{1>s_1>t_1>t_2>s_2>0}\right.\\
&\left.+\int\limits_{1>s_1>t_1>s_2>t_2>0}+\int\limits_{1>s_1>s_2>t_1>t_2>0}\right)\frac{dt_1dt_2ds_1ds_2}{t_1(1-t_2)s_1(1-s_2)}\\
=&2\zeta(2,2)+4\zeta(3,1).
\end{align*}
Such type products are called the shuffle products.

Combining the stuffle products and the shuffle products, one can obtain the (finite) double shuffle relations. For example, the relations
$$2\zeta(2,2)+\zeta(4)=\zeta(2)^2=2\zeta(2,2)+4\zeta(3,1)$$
belong to the system of the double shuffle relations. However, these are not all algebraic relations among multiple zeta values. For example, Euler's formula $\zeta(3)=\zeta(2,1)$ can not be deduced from the double shuffle relations. Hence we need to consider the divergent multiple zeta values and to take regularization processes. There are also two regularization processes, using the truncated multiple zeta value $\zeta_{N}(\mathbf{k})$
defined as a finite sum
$$\zeta_N(\mathbf{k})=\sum\limits_{N>m_1>\cdots>m_n>0}\frac{1}{m_1^{k_1}\cdots m_n^{k_n}}$$
with $N>0$ and the multiple polylogarithms
$\Li_{\mathbf{k}}(z)$ defined as an iterated integral
$$\Li_{\mathbf{k}}(z)=\int\limits_{z>t_1>\cdots>t_k>0}f_1(t_1)\cdots f_k(t_k)dt_1\cdots dt_k$$
with $|z|<1$, respectively. Here $\mathbf{k}=(k_1,\ldots,k_n)$ is a finite sequence of positive integers and the functions $f_i(t)$ are defined by \eqref{Eq:One-Forms}. Note that if $\mathbf{k}$ is admissible, we have the limits
$$\lim\limits_{N\rightarrow \infty}\zeta_N(\mathbf{k})=\zeta(\mathbf{k}),\quad \lim\limits_{z\rightarrow 1}\Li_{\mathbf{k}}(z)=\zeta(\mathbf{k}).$$
We have the asymptotic properties
\begin{align*}
\zeta_N(\mathbf{k})&=P_{\mathbf{k}}(\log
N+\gamma)+O(N^{-1}\log^{J}N),\quad N\rightarrow \infty\quad (\exists J>0),\\
\Li_{\mathbf{k}}(z)&=Q_{\mathbf{k}}(-\log(1-z))+O((1-z)\log^{J'}(1-z)),\quad z\rightarrow
1\quad (\exists J'>0),
\end{align*}
where $P_{\mathbf{k}}(T),Q_{\mathbf{k}}(T)\in\mathbb{R}[T]$ are
polynomials and $\gamma$ is Euler's
constant.  Let $\rho:\mathbb{R}[T]\rightarrow
\mathbb{R}[T]$ be the $\mathbb{R}$-linear map determined by
$$\rho(e^{Ts})=\exp\left(\sum\limits_{n=2}^\infty\frac{(-1)^n}{n}\zeta(n)s^n\right)e^{Ts},$$
where $s$ is a variable. Then for any finite sequence $\mathbf{k}$ of positive integers, we have
$$Q_{\mathbf{k}}(T)=\rho(P_{\mathbf{k}}(T)),$$
which deduces relations of multiple zeta values. Then the so called regularized double shuffle relations contain such type relations and the double shuffle relations. The regularized double shuffle relations can be expressed algebraically as K. Ihara, M. Kaneko and D. Zagier have done in \cite{Ihara-Kaneko-Zagier}.

Since the system of the regularized double shuffle relations is a conjectured full
relation system of multiple zeta values, all other algebraic relations should be consequences of this system. The purpose of this paper is to give as many as relations for multiple zeta values which can be derived from the regularized double shuffle relations. For example, for a system satisfying the regularized double shuffle relations, we reprove the weighted sum formula of L. Guo and B. Xie, give some evaluation formulas with even arguments and discuss the restricted sum formulas of M. E. Hoffman and their generalizations. Besides relations corresponding to multiple zeta values, we also provide in this system the formulas corresponding to multiple zeta-star values.

The paper is organized as follows. In Section \ref{Sec:DoubleShuffle}, we recall the algebraic setting of the regularized double shuffle relations. Besides descriptions corresponding to multiple zeta values, we also consider the descriptions corresponding to multiple zeta-star values. In Section \ref{Sec:Relations}, we give as many as the relations which can be deduced from the regularized double shuffle relations. Finally, in Section \ref{Sec:Remarks}, we briefly recall some other conjectured full relation systems of multiple zeta values.


\section{Algebraic setting of the regularized double shuffle relations}\label{Sec:DoubleShuffle}

In this section, we recall the algebraic setting of the regularized double shuffle relations from \cite{Ihara-Kaneko-Zagier} for multiple zeta values and from \cite{Muneta} for multiple zeta-star values. Some new materials are supplied.

Let $\mathbb{K}$ be a field of characteristic zero and $R$ be a commutative $\mathbb{K}$-algebra with unitary.


\subsection{Double shuffle relations}

Let $A=\{x,y\}$ be an alphabet containing two noncommutative  letters. Denote by $A^{\ast}$ the set of all words generated by the letters in $A$, which contains the empty word $1$. The subset of the nonempty words of $A^{\ast}$ is denoted by $A^{+}$. Hence $A^{+}=A^{\ast}-\{1\}$. Let $\mathfrak{h}=\mathbb{K}\langle A\rangle$ be the $\mathbb{K}$-algebra of noncommutative polynomials on $x$ and $y$ with coefficients in $\mathbb{K}$. Then $A^{\ast}$ is a basis of the $\mathbb{K}$-vector space $\mathfrak{h}$. We define two subalgebras
$$\mathfrak{h}^1=\mathbb{K}+\mathfrak{h}y,\quad \mathfrak{h}^0=\mathbb{K}+x\mathfrak{h}y,$$
which have the $\mathbb{K}$-basis
$$z_{k_1}\cdots z_{k_n},\quad n\geqslant 0,k_1,\ldots,k_n\geqslant 1$$
and
$$z_{k_1}\cdots z_{k_n},\quad n\geqslant 0,k_1\geqslant 2,k_2,\ldots,k_n\geqslant 1,$$
respectively, where $z_k=x^{k-1}y$ for any positive integer $k$.

We define a new product $\shuffle$ on $\mathfrak{h}$ by the $\mathbb{K}$-bilinearities and the axioms
\begin{itemize}
  \item [(S1)] $w\shuffle 1=1\shuffle w=w$, \quad ($\forall w\in A^{\ast}$);
  \item [(S2)] $aw\shuffle bv=a(w\shuffle bv)+b(aw\shuffle v)$,\quad ($\forall a,b\in A,\forall w,v\in A^{\ast}$).
\end{itemize}
The product $\shuffle$ is commutative and associative, and is called  the shuffle product. Under this new product, $\mathfrak{h}$ becomes a commutative algebra, and $\mathfrak{h}^1$ and $\mathfrak{h}^0$ are still subalgebras. We call these commutative algebras shuffle algebras, and denote them by $\mathfrak{h}_{\shuffle}$, $\mathfrak{h}^1_{\shuffle}$ and $\mathfrak{h}^0_{\shuffle}$, respectively.

Another product $\ast$ can be defined on $\mathfrak{h}^1$ by the $\mathbb{K}$-bilinearities and the axioms
\begin{itemize}
  \item [(H1)] $w\ast 1=1\ast w=w$, \quad ($\forall w\in A^{\ast}\cap \mathfrak{h}^1$);
  \item [(H2)] $z_kw\ast z_lv=z_k(w\ast z_lv)+z_l(z_kw\ast v)+z_{k+l}(w\ast v)$,\; ($\forall k,l\in \mathbb{N},\forall w,v\in A^{\ast}\cap\mathfrak{h}^1$).
\end{itemize}
Here $\mathbb{N}$ is the set of all positive integers. The product $\ast$ is commutative and associative, and is called  the stuffle product (also called harmonic product). Under this new product, $\mathfrak{h}^1$ becomes a commutative algebra, and $\mathfrak{h}^0$ is still a subalgebra. We call these commutative algebras stuffle algebras, and denote them by $\mathfrak{h}^1_{\ast}$ and $\mathfrak{h}^0_{\ast}$, respectively.

Let $\sigma$ be the automorphism of the noncommutative algebra $\mathfrak{h}$ determined by
$$\sigma(x)=x,\quad \sigma(y)=x+y.$$
We know that the inverse map $\sigma^{-1}$ is determined by
$$\sigma^{-1}(x)=x,\quad \sigma^{-1}(y)=-x+y.$$
Note that $\sigma$ and $\sigma^{-1}$ are also automorphisms of the shuffle algebra $\mathfrak{h}_{\shuffle}$.
We define a $\mathbb{K}$-linear map $S:\mathfrak{h}\rightarrow \mathfrak{h}$, such that $S(1)=1$ and
$$S(wa)=\sigma(w)a,\quad (\forall a\in A,\forall w\in A^{\ast}).$$
Then it is easy to see that $S$ is invertible, and the inverse linear map $S^{-1}$ is given by $S^{-1}(1)=1$ and
$$S^{-1}(wa)=\sigma^{-1}(w)a,\quad (\forall a\in A,\forall w\in A^{\ast}).$$
We have
$$S(\mathfrak{h}^1)=\mathfrak{h}^1,\quad S(\mathfrak{h}^0)=\mathfrak{h}^0,$$
which imply that we have invertible linear maps $S|_{\mathfrak{h}^1}$ and $S|_{\mathfrak{h}^0}$.
Note that the map $S$ defined in \cite{Muneta} is in fact $S|_{\mathfrak{h}^1}$ here (See also \cite{Li-Qin}).

We define the shuffle product $\sshuffle$ associated to multiple zeta-star values on $\mathfrak{h}$ by the $\mathbb{K}$-bilinearities and the axioms
\begin{itemize}
  \item [(SS1)] $w\sshuffle 1=1\sshuffle w=w$, \quad ($\forall w\in A^{\ast}$);
  \item [(SS2)] $aw\sshuffle bv=a(w\sshuffle bv)+b(aw\sshuffle v)-\delta(w)\rho(a)bv-\delta(v)\rho(b)aw$,

  ($\forall a,b\in A,\forall w,v\in A^{\ast}$).
\end{itemize}
where $\delta:A^{\ast}\rightarrow \{0,1\}$ is a map defined by
$$\delta(w)=\begin{cases}
1 & \text{if\;}w=1,\\
0 & \text{if\;}w\neq 1,
\end{cases}$$
and $\rho:A\longrightarrow \mathfrak{h}$ is a map defined by
$$\rho(x)=0,\quad \rho(y)=x.$$
Note that in \cite{Muneta}, S. Muneta used a different map instead of $\rho$. By \cite{Li-Qin}, for any $w,v\in\mathfrak{h}$, we have
$$S(w\sshuffle v)=S(w)\shuffle S(v),\quad S^{-1}(w\shuffle v)=S^{-1}(w)\sshuffle S^{-1}(v),$$
which imply that the product $\sshuffle$ is also commutative and associative. Hence we get the commutative algebra $\mathfrak{h}_{\sshuffle}$ and its subalgebras $\mathfrak{h}_{\sshuffle}^1$ and $\mathfrak{h}^0_{\sshuffle}$.

As in \cite{Muneta}, the stuffle product $\sast$ associated to multiple zeta-star values on $\mathfrak{h}^1$ is defined by the $\mathbb{K}$-bilinearities and the axioms
\begin{itemize}
  \item [(SH1)] $w\sast 1=1\sast w=w$, \quad ($\forall w\in A^{\ast}\cap \mathfrak{h}^1$);
  \item [(SH2)] $z_kw\sast z_lv=z_k(w\sast z_lv)+z_l(z_kw\sast v)-z_{k+l}(w\sast v)$,\; ($\forall k,l\in \mathbb{N},\forall w,v\in A^{\ast}\cap\mathfrak{h}^1$).
\end{itemize}
By \cite{Muneta}, for any $w,v\in\mathfrak{h}^1$, we have
$$S(w\sast v)=S(w)\ast S(v),\quad S^{-1}(w\ast v)=S^{-1}(w)\sast S^{-1}(v),$$
which imply that the product $\sast$ is also commutative and associative. Hence we get the commutative algebra $\mathfrak{h}^1_{\sast}$ and its subalgebra $\mathfrak{h}^0_{\sast}$.

Now for any $\mathbb{K}$-linear map $Z_R:\mathfrak{h}^0\rightarrow R$, we define
$$Z_R^{\star}=Z_R\circ S:\mathfrak{h}^0\rightarrow R,$$
which is also a $\mathbb{K}$-linear map. Then for any admissible index $\mathbf{k}=(k_1,\ldots,k_n)$, we define the multiple zeta value $\zeta_R(\mathbf{k})$ and
the multiple zeta-star value $\zeta^{\star}_R(\mathbf{k})$ associated with the map $Z_R$ respectively by
$$\zeta_R(\mathbf{k})=Z_R(z_{k_1}\cdots z_{k_n}),\qquad \zeta^{\star}_R(\mathbf{k})=Z_R^{\star}(z_{k_1}\cdots z_{k_n}).$$

The following lemma is easy to prove.

\begin{lem}\label{Lem:FiniteDoubleShuffle}
For any $\mathbb{K}$-linear map $Z_R:\mathfrak{h}^0\rightarrow R$, the following two properties are equivalent
\begin{itemize}
  \item [(i)] $Z_R:\mathfrak{h}^0_{\shuffle}\rightarrow R$ and $Z_R:\mathfrak{h}^0_{\ast}\rightarrow R$ are algebra homomorphisms, that is for any $w,v\in\mathfrak{h}^0$, we have
      $$Z_R(w\shuffle v)=Z_R(w)Z_R(v)=Z_R(w\ast v);$$
  \item [(ii)] $Z_R^{\star}:\mathfrak{h}^0_{\sshuffle}\rightarrow R$ and $Z_R^{\star}:\mathfrak{h}^0_{\sast}\rightarrow R$ are algebra homomorphisms, that is for any $w,v\in\mathfrak{h}^0$, we have
      $$Z_R^{\star}(w\sshuffle v)=Z_R^{\star}(w)Z_R^{\star}(v)=Z_R^{\star}(w\sast v).$$
\end{itemize}
\end{lem}

If a $\mathbb{K}$-linear map $Z_R:\mathfrak{h}^0\rightarrow R$ satisfies the equivalent properties in Lemma \ref{Lem:FiniteDoubleShuffle}, we call the pair $(R,Z_R)$ or simply the map $Z_R$ satisfies the (finite) double shuffle relations. Below, if $Z_R$ satisfies the double shuffle relations, we always assume that $Z_R(1)=1$. Hence if $\mathbf{k}$ is an empty index, we set
$$\zeta_R(\mathbf{k})=\zeta_R^{\star}(\mathbf{k})=Z_R(1)=1.$$


\subsection{Regularized double shuffle relations}

By \cite{Reutenauer}, the algebra $\mathfrak{h}^1_{\shuffle}$ is a polynomial algebra $\mathfrak{h}^1_{\shuffle}=\mathfrak{h}^0_{\shuffle}[y]$. Then we can define the regularization map $\reg_{\shuffle}:\mathfrak{h}^1_{\shuffle}\rightarrow \mathfrak{h}^0_{\shuffle}$, such that for any $$w=\sum\limits_{i\geqslant 0}w_i\shuffle y^{\shuffle i}\in\mathfrak{h}^1$$
with $w_i\in\mathfrak{h}^0$, we have $\reg_{\shuffle}(w)=w_0$. The map $\reg_{\shuffle}$ is an algebra homomorphism. Similarly, by \cite{Hoffman1997}, the algebra $\mathfrak{h}^1_{\ast}$ is also a polynomial algebra $\mathfrak{h}^1_{\ast}=\mathfrak{h}^0_{\ast}[y]$. Then we can define the regularization map $\reg_{\ast}:\mathfrak{h}^1_{\ast}\rightarrow \mathfrak{h}^0_{\ast}$, which is an algebra homomorphism.

For any $w\in\mathfrak{h}^1$, we can write $S(w)\in\mathfrak{h}^1$ as
$$S(w)=\sum\limits_{i\geqslant 0}w_i\shuffle y^{\shuffle i}$$
with $w_i\in\mathfrak{h}^0$. Applying the map $S^{-1}$, we get
$$w=\sum\limits_{i\geqslant 0}S^{-1}(w_i)\sshuffle y^{\sshuffle i}.$$
Here $S^{-1}(w_i)\in \mathfrak{h}^0$. Furthermore, it is easy to see that such decomposition is unique. Hence $\mathfrak{h}^{1}_{\sshuffle}=\mathfrak{h}^{0}_{\sshuffle}[y]$ is also a polynomial algebra, and one can define the regularization map $\reg_{\sshuffle}:\mathfrak{h}^1_{\sshuffle}\rightarrow \mathfrak{h}^0_{\sshuffle}$, which is an algebra homomorphism. Similarly, $\mathfrak{h}^1_{\sast}=\mathfrak{h}^0_{\sast}[y]$, and one can define an algebra homomorphism $\reg_{\sast}:\mathfrak{h}^1_{\sast}\rightarrow\mathfrak{h}^0_{\sast}$. Note that we have
$$\reg_{\sshuffle}=S^{-1}\circ \reg_{\shuffle}\circ S,\qquad \reg_{\sast}=S^{-1}\circ \reg_{\ast}\circ S.$$
Here and below, both $S|_{\mathfrak{h}^1}$ and $S|_{\mathfrak{h}^0}$ are denoted by $S$.

Now assume that the $\mathbb{K}$-linear map $Z_{R}:\mathfrak{h}^{0}\rightarrow R$ satisfies the double shuffle relations. Then there exist unique algebra homomorphisms
\begin{align*}
&Z_{R}^{\shuffle}:\mathfrak{h}^{1}_{\shuffle}\rightarrow R[T],\quad  Z_{R}^{\ast}:\mathfrak{h}^{1}_{\ast}\rightarrow R[T],\\
&Z_{R}^{\sshuffle}:\mathfrak{h}^{1}_{\sshuffle}\rightarrow R[T],\quad Z_{R}^{\sast}:\mathfrak{h}^{1}_{\sast}\rightarrow R[T],
\end{align*}
such that
\begin{align*}
&\left.Z_{R}^{\shuffle}\right|_{\mathfrak{h}^{0}}=Z_{R}=\left.Z_{R}^{\ast}\right|_{\mathfrak{h}^{0}},\qquad
Z_{R}^{\shuffle}(y)=T=Z_{R}^{\ast}(y),\\
&\left.Z_{R}^{\sshuffle}\right|_{\mathfrak{h}^{0}}=Z_{R}^{\star}=\left.Z_{R}^{\sast}\right|_{\mathfrak{h}^{0}},\qquad
Z_{R}^{\sshuffle}(y)=T=Z_{R}^{\sast}(y).
\end{align*}
More precisely, assume that $w\in\mathfrak{h}^1$ is decomposed as
$$w=\sum\limits_{i\geqslant 0}w_i\shuffle y^{\shuffle i}=\sum\limits_{j\geqslant 0}\widetilde{w}_j\ast y^{\ast j},\quad (w_i,\widetilde{w}_j\in\mathfrak{h}^0),$$
then we have
$$Z_R^{\shuffle}(w)=\sum\limits_{i\geqslant 0}Z_R(w_i)T^i,\quad Z_R^{\ast}(w)=\sum\limits_{j\geqslant 0}Z_R(\widetilde{w}_j)T^j;$$
and if $v\in\mathfrak{h}^1$ is decomposed as
$$v=\sum\limits_{i\geqslant 0}v_i\sshuffle y^{\sshuffle i}=\sum\limits_{j\geqslant 0}\widetilde{v}_j\sast y^{\sast j},\quad (v_i,\widetilde{v}_j\in\mathfrak{h}^0),$$
then we have
$$Z_R^{\sshuffle}(v)=\sum\limits_{i\geqslant 0}Z_R^{\star}(v_i)T^i,\quad Z_R^{\sast}(v)=\sum\limits_{j\geqslant 0}Z_R^{\star}(\widetilde{v}_j)T^j.$$
Finally by the definitions and the properties of the map $S$, we get
$$\begin{array}{l@{\qquad\qquad}l}
Z_R^{\sshuffle}=Z_R^{\shuffle}\circ S, & Z_R^{\sast}=Z_R^{\ast}\circ S; \\
Z_R^{\shuffle}|_{T=0}=Z_R\circ \reg_{\shuffle}, & Z_R^{\ast}|_{T=0}=Z_R\circ \reg_{\ast};\\
Z_R^{\sshuffle}|_{T=0}=Z_R^{\star}\circ \reg_{\sshuffle}, & Z_R^{\sast}|_{T=0}=Z_R^{\star}\circ \reg_{\sast}.
\end{array}$$

For any $\mathbb{K}$-linear map $Z_R:\mathfrak{h}^0\rightarrow R$, let
$$\Gamma_{R}(s)=\exp\left(\sum\limits_{n=2}^{\infty}\frac{(-1)^{n}}{n}\zeta_R(n)s^n\right)\in R[[s]],$$
which is called the gamma series associated to the map $Z_R$.
Let $\rho_{R}:R[T]\rightarrow R[T]$ be the $R$-module homomorphism defined by
$$\rho_R(e^{Ts})=\Gamma_{R}(s)e^{Ts}.$$
Note that $\rho_R$ is invertible, and the inverse $\rho_R^{-1}$ is determined by
$$\rho_R^{-1}(e^{Ts})=\Gamma_{R}(s)^{-1}e^{Ts}.$$

In \cite{Ihara-Kaneko-Zagier}, K. Ihara, M. Kaneko and D. Zagier showed that one can use derivations to describe the regularized double shuffle relations for multiple zeta values. Recall that a $\mathbb{K}$-linear map $\mathcal{D}:\mathfrak{h}\rightarrow \mathfrak{h}$ is a derivation if it satisfies
$$\mathcal{D}(wv)=\mathcal{D}(w)v+w\mathcal{D}(v),\quad (\forall w,v\in\mathfrak{h}).$$
Then we see that $\mathcal{D}(1)=0$. To deal with the multiple zeta-star values, we introduce the following definitions.

\begin{dfn}
Assume that $\mathcal{D}$ and $\mathcal{L}$ are $\mathbb{K}$-linear maps on $\mathfrak{h}$ with $\mathcal{D}(1)=0$ and $\mathcal{L}$ invertible. If
$$\mathcal{D}(wv)=(\mathcal{L}^{-1}\circ \mathcal{D}\circ \mathcal{L})(w)v+w\mathcal{D}(v),\quad (\forall w,v\in A^{+}),$$
we call $\mathcal{D}$ is a left $\mathcal{L}$-derivation on $\mathfrak{h}$. Similarly, if
$$\mathcal{D}(wv)=\mathcal{D}(w)v+w(\mathcal{L}^{-1}\circ \mathcal{D}\circ \mathcal{L})(v),\quad (\forall w,v\in A^{+}),$$
we call $\mathcal{D}$ is a right $\mathcal{L}$-derivation on $\mathfrak{h}$.
\end{dfn}

Note that if $\mathcal{D}$ is a left $\mathcal{L}$-derivation, then for any $a_1,\ldots, a_n\in A$, we have
\begin{align*}
\mathcal{D}(a_1\cdots a_n)=&\sum\limits_{i=1}^{n-1}a_1\cdots a_{i-1}(\mathcal{L}^{-1}\circ \mathcal{D}\circ \mathcal{L})(a_i)a_{i+1}\cdots a_n\\
&+a_1\cdots a_{n-1}\mathcal{D}(a_n).
\end{align*}

We define a $\mathbb{K}$-linear map $\widetilde{S}:\mathfrak{h}\rightarrow \mathfrak{h}$ such that $\widetilde{S}(1)=1$ and
$$\widetilde{S}(wa)=w\sigma(a),\quad (\forall a\in A,\forall w\in A^{\ast}).$$
Then it is easy to check that $\widetilde{S}=S^{-1}\circ \sigma$, which implies that $\widetilde{S}$ is invertible with inverse $\widetilde{S}^{-1}=\sigma^{-1}\circ S$, and
$$\widetilde{S}^{-1}(wa)=w\sigma^{-1}(a),\quad (\forall a\in A,\forall w\in A^{\ast}).$$
Then the following result holds.

\begin{lem}\label{Lem:Left-Derivation}
Assume $\mathcal{D}:\mathfrak{h}\rightarrow\mathfrak{h}$ is a $\mathbb{K}$-linear map, such that $\mathcal{D}(w)\not\in \mathbb{K}$ for any $w\in A^{+}$. Set $\mathcal{D}^{\star}=S^{-1}\circ \mathcal{D}\circ S$. Then $\mathcal{D}$ is a derivation if and only if $\mathcal{D}^{\star}$ is a left $\widetilde{S}$-derivation.
\end{lem}

\proof
Assume that $\mathcal{D}$ is a derivation on $\mathfrak{h}$. For any $w,v\in A^{+}$, we have
$$\mathcal{D}^{\star}(wv)=S^{-1}(\mathcal{D}(\sigma(w)S(v)))=S^{-1}(\mathcal{D}(\sigma(w))S(v)+\sigma(w)\mathcal{D}(S(v))).$$
Since $S(v),\mathcal{D}(S(v))\not\in\mathbb{K}$, we have
\begin{align*}
&\mathcal{D}^{\star}(wv)=\sigma^{-1}(\mathcal{D}(\sigma(w)))v+w S^{-1}(\mathcal{D}(S(v)))\\
=&(\sigma^{-1}\circ S\circ \mathcal{D}^{\star}\circ S^{-1}\circ \sigma)(w)v+w\mathcal{D}^{\star}(v)\\
=&(\widetilde{S}^{-1}\circ\mathcal{D}^{\star}\circ \widetilde{S})(w)v+w\mathcal{D}^{\star}(v).
\end{align*}
Hence $\mathcal{D}^{\star}$ is a left $\widetilde{S}$-derivation. one can prove the oppositive direction similarly.
\qed

For any $\mathbb{K}$-linear map $\mathcal{D}$ on $\mathfrak{h}$, let $\overline{\mathcal{D}}=\tau\circ\mathcal{D}\circ\tau$. Here $\tau$ is the antiautomorphism of the noncommutative algebra $\mathfrak{h}$ determined by
$$\tau(x)=y,\qquad \tau(y)=x.$$
Then $\overline{\mathcal{D}}$ is also $\mathbb{K}$-linear. Note that $\tau$ is an automorphism of the shuffle algebra $\mathfrak{h}_{\shuffle}$. We have the following interesting results, although we will not use them below.

\begin{lem}
Assume that $\mathcal{D}$ and $\mathcal{L}$ are $\mathbb{K}$-linear maps on $\mathfrak{h}$ with $\mathcal{L}$ invertible. Then
\begin{itemize}
  \item [(1)] $\mathcal{D}$ is a derivation if and only if $\overline{\mathcal{D}}$ is a derivation;
  \item [(2)] $\mathcal{D}$ is a left $\mathcal{L}$-derivation if and only if $\overline{\mathcal{D}}$ is a right $\overline{\mathcal{L}}$-derivation.
\end{itemize}
\end{lem}

\proof (1) Assume that $\mathcal{D}$ is a derivation on $\mathfrak{h}$. For any $w,v\in\mathfrak{h}$, we have
\begin{align*}
&\overline{\mathcal{D}}(wv)=\tau(\mathcal{D}(\tau(v)\tau(w)))=\tau(\mathcal{D}(\tau(v))\tau(w)+\tau(v)\mathcal{D}(\tau(w)))\\
=&w\tau(\mathcal{D}(\tau(v)))+\tau(\mathcal{D}(\tau(w)))v=w\overline{\mathcal{D}}(v)+\overline{\mathcal{D}}(w)v.
\end{align*}
Hence $\overline{\mathcal{D}}$ is a derivation. On the contrary, using $\mathcal{D}=\overline{\overline{\mathcal{D}}}$ we get the oppositive direction.

\noindent (2) Assume that $\mathcal{D}$ is a left $\mathcal{L}$-derivation on $\mathfrak{h}$. For any $w,v\in A^{+}$, we have
\begin{align*}
&\overline{\mathcal{D}}(wv)=\tau(\mathcal{D}(\tau(v)\tau(w)))=\tau((\mathcal{L}^{-1}\circ\mathcal{D}\circ \mathcal{L})(\tau(v))\tau(w)+\tau(v)\mathcal{D}(\tau(w)))\\
=&w(\tau\circ\mathcal{L}^{-1}\circ \mathcal{D}\circ\mathcal{L}\circ \tau)(v)+\overline{\mathcal{D}}(w)v\\
=&w(\overline{\mathcal{L}}^{-1}\circ\overline{\mathcal{D}}\circ\overline{\mathcal{L}})(v)+\overline{\mathcal{D}}(w)v.
\end{align*}
Hence $\overline{\mathcal{D}}$ is a right $\overline{\mathcal{L}}$-derivation.
\qed

Now we define the derivations which describe the regularized double shuffle relations. For any positive integer $n$, the derivation $\partial_n$ on $\mathfrak{h}$ is determined by
$$\partial_n(x)=x(x+y)^{n-1}y,\quad \partial_n(y)=-x(x+y)^{n-1}y.$$
Set $\partial_n^{\star}=S^{-1}\circ \partial_n \circ S$. Then $\partial_n^{\star}$ is a left $\widetilde{S}$-derivation on $\mathfrak{h}$, and satisfies
$$\partial_n^{\star}(x)=xy^n,\quad \partial_n^{\star}(y)=-xy^n.$$

We need some other maps, which have a close relation with Ohno's relations (\cite{Ohno}). For any nonnegative integer $m$, the $\mathbb{K}$-linear map $\sigma_m:\mathfrak{h}^1\rightarrow\mathfrak{h}^1$ is defined by $\sigma_m(1)=1$, and
$$\sigma_m(z_{k_1}\cdots z_{k_n})=\sum\limits_{\varepsilon_1+\cdots+\varepsilon_n=m\atop \varepsilon_i\geqslant 0}z_{k_1+\varepsilon_1}\cdots z_{k_n+\varepsilon_n},\quad (\forall n,k_1,\ldots,k_n\in\mathbb{N}).$$
Set $\sigma_m^{\star}=S^{-1}\circ \sigma_m\circ S$ and $\overline{\sigma_m}^{\star}=S^{-1}\circ \overline{\sigma_m}\circ S$. Finally, we can state the equivalent properties about the regularized double shuffle relations.

\begin{thm}\label{Thm:RegDoubleShuffle}
Assume that the pair $(R,Z_R)$ satisfies the double shuffle relations. Then the following properties are equivalent:
\begin{itemize}
    \item[(i)] $(Z_{R}^{\shuffle}-\rho_{R}\circ Z_{R}^{\ast})(w_1)=0$ for all $w_1\in\mathfrak{h}^{1}$;
    \item[(a)] $(Z_{R}^{\sshuffle}-\rho_{R}\circ Z_{R}^{\sast})(w_1)=0$ for all $w_1\in\mathfrak{h}^{1}$;
    \item[(ii)] $(Z_{R}^{\shuffle}-\rho_{R}\circ Z_{R}^{\ast})(w_1)|_{T=0}=0$ for all $w_1\in\mathfrak{h}^{1}$;
     \item[(b)] $(Z_{R}^{\sshuffle}-\rho_{R}\circ Z_{R}^{\sast})(w_1)|_{T=0}=0$ for all $w_1\in\mathfrak{h}^{1}$;
    \item[(iii)] $Z_{R}^{\shuffle}(w_1\shuffle w_0-w_1\ast w_0)=0$ for all $w_1\in\mathfrak{h}^{1}$ and all $w_0\in\mathfrak{h}^{0}$;
    \item[(c)] $Z_{R}^{\sshuffle}(w_1\sshuffle w_0-w_1\sast w_0)=0$ for all $w_1\in\mathfrak{h}^{1}$ and all $w_0\in\mathfrak{h}^{0}$;
    \item[(iii$'$)] $Z_{R}^{\ast}(w_1\shuffle w_0-w_1\ast w_0)=0$ for all $w_1\in\mathfrak{h}^{1}$ and all $w_0\in\mathfrak{h}^{0}$;
    \item[(c$'$)] $Z_{R}^{\sast}(w_1\sshuffle w_0-w_1\sast w_0)=0$ for all $w_1\in\mathfrak{h}^{1}$ and all $w_0\in\mathfrak{h}^{0}$;
    \item[(iv)] $Z_{R}(\reg_{\shuffle}(w_1\shuffle w_0-w_1\ast w_0))=0$ for all $w_1\in\mathfrak{h}^{1}$ and all $w_0\in\mathfrak{h}^{0}$;
    \item[(d)] $Z_{R}^{\star}(\reg_{\sshuffle}(w_1\sshuffle w_0-w_1\sast w_0))=0$ for all $w_1\in\mathfrak{h}^{1}$ and all $w_0\in\mathfrak{h}^{0}$;
    \item[(iv$'$)] $Z_{R}(\reg_{\ast}(w_1\shuffle w_0-w_1\ast w_0))=0$ for all $w_1\in\mathfrak{h}^{1}$ and all $w_0\in\mathfrak{h}^{0}$;
    \item[(d$'$)] $Z_{R}^{\star}(\reg_{\sast}(w_1\sshuffle w_0-w_1\sast w_0))=0$ for all $w_1\in\mathfrak{h}^{1}$ and all $w_0\in\mathfrak{h}^{0}$;
    \item[(v)] $Z_{R}(\reg_{\shuffle}(y^n\ast w_0))=0$ for all positive integer $n$ and all $w_0\in \mathfrak{h}^{0}$;
    \item[(e)] $Z_{R}^{\star}(\reg_{\sshuffle}(y^n\sshuffle w_0-y^n\sast w_0))=0$ for all positive integer $n$ and all $w_0\in \mathfrak{h}^{0}$;
    \item[(v$'$)] $Z_{R}(\reg_{\ast}(y^n\shuffle w_0-y^n\ast w_0))=0$ for all positive integer $n$ and all $w_0\in \mathfrak{h}^{0}$;
    \item[(e$'$)] $Z_{R}^{\star}(\reg_{\sast}(y^n\sshuffle w_0-y^n\sast w_0))=0$ for all positive integer $n$ and all $w_0\in \mathfrak{h}^{0}$;
    \item[(vi)] $Z_{R}(\partial_n(w_0))=0$ for all positive integer $n$ and all $w_0\in \mathfrak{h}^{0}$;
    \item[(f)] $Z_{R}^{\star}(\partial_n^{\star}(w_0))=0$ for all positive integer $n$ and all $w_0\in \mathfrak{h}^{0}$;
    \item[(vii)] $Z_{R}((\sigma_{m}-\overline{\sigma_{m}})(w_0))=0$ for all nonnegative integer $m$ and all $w_0\in \mathfrak{h}^{0}$;
    \item[(g)] $Z_{R}^{\star}((\sigma_{m}^{\star}-\overline{\sigma_{m}}^{\star})(w_0))=0$ for all nonnegative integer $m$ and all $w_0\in \mathfrak{h}^{0}$.
\end{itemize}
\end{thm}

\begin{rem}
The items (i)-(vii) and (iii$'$)-(v$'$) in Theorem \ref{Thm:RegDoubleShuffle} were listed in \cite[Theorem 2 and Theorem 3]{Ihara-Kaneko-Zagier}. Here we add the equivalent properties corresponding to the multiple zeta-star values, which are marked by (a), (b), and so on.
\end{rem}

If the $\mathbb{K}$-linear map $Z_R:\mathfrak{h}^0\longrightarrow R$ satisfies the double shuffle relations, and satisfies the equivalent properties in Theorem \ref{Thm:RegDoubleShuffle}, we call the pair $(R,Z_R)$ or simply the map $Z_R$ satisfies the regularized double shuffle relations.
Hence by \cite{Ihara-Kaneko-Zagier}, if we take $\mathbb{K}=\mathbb{Q}$, $R=\mathbb{R}$, and define the $\mathbb{Q}$-linear map $Z_{R}=Z:\mathfrak{h}^{0}\rightarrow \mathbb{R}$, such that
$$Z(z_{k_1}z_{k_2}\cdots z_{k_n})=\zeta(k_1,k_2,\cdots,k_n),\quad (k_1\geqslant 2,k_2,\ldots,k_n\geqslant 1),$$
then the pair $(\mathbb{R},Z)$ satisfies the regularized double shuffle relations.

To prove Theorem \ref{Thm:RegDoubleShuffle}, we need the following lemmas.

\begin{lem}\label{Lem:ym-ymw}
Assume that the pair $(R,Z_R)$ satisfies the double shuffle relations and $n$ is a positive integer. If
$$Z_R(\reg_{\odot}(y^m\shuffle w_0-y^m\ast w_0))=0$$
for any integer $m$ with $1\leqslant m\leqslant n$ and any $w_0\in\mathfrak{h}^0$, then
$$Z_R(\reg_{\odot}(y^mw_1\shuffle w_2-y^mw_1\ast w_2))=0$$
for any integer $m$ with $1\leqslant m\leqslant n$ and any $w_1,w_2\in\mathfrak{h}^0$, where $\odot=\shuffle$ or $\ast$.
\end{lem}

\proof
We prove the case of $\odot=\shuffle$, and one can prove the case of  $\odot=\ast$ similarly. We process by induction on $n$. For $n=1$, note that for any $w_1,w_2\in\mathfrak{h}^0$, we have
$$yw_1-y\shuffle w_1,yw_1-y\ast w_1\in\mathfrak{h}^0,$$
and
\begin{align*}
&yw_1\shuffle w_2-yw_1\ast w_2=(yw_1-y\shuffle w_1)\shuffle w_2-(yw_1-y\ast w_1)\ast w_2\\
&\qquad +y\shuffle (w_1\ast w_2)-y\ast (w_1\ast w_2)+y\shuffle (w_1\shuffle w_2-w_1\ast w_2).
\end{align*}
Hence we get
\begin{align*}
&Z_R(\reg_{\shuffle}(yw_1\shuffle w_2-yw_1\ast w_2))\\
=&Z_R(yw_1-y\shuffle w_1)Z_R(w_2)-Z_R(yw_1-y\ast w_1)Z_R(w_2)\\
&\quad+Z_R(\reg_{\shuffle}(y\shuffle (w_1\ast w_2)-y\ast (w_1\ast w_2)))\\
&\quad +Z_R(\reg_{\shuffle}(y))Z_R(w_1\shuffle w_2-w_1\ast w_2)\\
=&Z_R(y\ast w_1-y\shuffle w_1)Z_R(w_2)=0.
\end{align*}

Now assume that $n>1$. For $w_1,w_2\in\mathfrak{h}^0$, there exist $u_i\in\mathfrak{h}^0$ such that
\begin{align*}
&y^nw_1\shuffle w_2-y^nw_1\ast w_2\\
=&(y^nw_1-y^n\ast w_1)\shuffle w_2-(y^nw_1-y^n\ast w_1)\ast w_2\\
&\quad-(y^n\shuffle w_1-y^n\ast w_1)\shuffle w_2\\
&\quad +y^n\shuffle(w_1\ast w_2)-y^n\ast(w_1\ast w_2)\\
&\quad +y^n\shuffle(w_1\shuffle w_2-w_1\ast w_2)\\
=&\sum\limits_{i=0}^{n-1}(y^iu_i\shuffle w_2-y^iu_i\ast w_2)-(y^n\shuffle w_1-y^n\ast w_1)\shuffle w_2\\
&+y^n\shuffle(w_1\ast w_2)-y^n\ast(w_1\ast w_2)+y^n\shuffle(w_1\shuffle w_2-w_1\ast w_2).
\end{align*}
Hence using the induction hypothesis, we easily get
$$Z_R(\reg_{\shuffle}(y^nw_1\shuffle w_2-y^nw_1\ast w_2))=0,$$
which finishes the proof.
\qed

\begin{lem}\label{Lem:Sym-ym}
Assume that the pair $(R,Z_R)$ satisfies the double shuffle relations and $n$ is a positive integer. If
$$Z_R(\reg_{\odot}(S(y^m)\shuffle w_0-S(y^m)\ast w_0))=0$$
for any integer $m$ with $1\leqslant m\leqslant n$ and any $w_0\in\mathfrak{h}^0$, then
$$Z_R(\reg_{\odot}(y^m\shuffle w_1-y^m\ast w_1))=0$$
for any integer $m$ with $1\leqslant m\leqslant n$ and any $w_1\in\mathfrak{h}^0$, where $\odot=\shuffle$ or $\ast$.
\end{lem}

\proof
We prove the case of $\odot=\shuffle$, and one can prove the case of  $\odot=\ast$ similarly. We process by induction on $n$. For $n=1$, the result follows from the fact that $S(y)=y$. Now assume that $n>1$. Since
$$S(y^n)=(x+y)^{n-1}y=\sum\limits_{i=0}^{n-2}y^iv_i+y^n$$
with $v_i\in\mathfrak{h}^0$, we have
$$S(y^n)\shuffle w_1-S(y^n)\ast w_1=\sum\limits_{i=0}^{n-2}(y^iv_i\shuffle w_1-y^iv_i\ast w_1)+(y^n\shuffle w_1-y^n\ast w_1).$$
Using the induction hypothesis, we have
$$Z_R(\reg_{\shuffle}(y^m\shuffle w_2-y^m\ast w_2))=0$$
for any integer $m$ with $1\leqslant m\leqslant n-1$ and any $w_2\in\mathfrak{h}^0$. Then by Lemma \ref{Lem:ym-ymw}, we get
$$Z_R(\reg_{\shuffle}(y^iv_i\shuffle w_1-y^iv_i\ast w_1))=0$$
for any $i$ with the condition $0\leqslant i\leqslant n-2$. Therefore we have
$$Z_R(\reg_{\shuffle}(y^n\shuffle w_1-y^n\ast w_1))=Z_R(\reg_{\shuffle}(S(y^n)\shuffle w_1-S(y^n)\ast w_1))=0,$$
which finishes the proof.
\qed

Now we give a proof of Theorem \ref{Thm:RegDoubleShuffle}.

\noindent {\bf Proof of Theorem \ref{Thm:RegDoubleShuffle}.} One can easily show the implications
$$\begin{array}{ll}
\text{(i)} \Longleftrightarrow \text{(a)},\quad \text{(ii)} \Longleftrightarrow \text{(b)}, & \text{(iii)} \Longleftrightarrow \text{(c)},\quad \text{(iii$'$\text)} \Longleftrightarrow \text{(c$'$)},\\
\text{(iv)} \Longleftrightarrow \text{(d)} \Longrightarrow \text{(e)},& \text{(iv$'$)} \Longleftrightarrow \text{(d$'$)} \Longrightarrow \text{(e$'$)},\\
\text{(vi)} \Longleftrightarrow \text{(f)}, & \text{(vii)} \Longleftrightarrow \text{(g)}.
\end{array}$$
Hence with the help of \cite[Theorem 2 and Theorem 3]{Ihara-Kaneko-Zagier}, to finish the proof, it is enough to show the implications
$$\text{(e)} \Longrightarrow \text{(v)},\quad \text{(e$'$)} \Longrightarrow \text{(v$'$)}.$$
Here we give a proof of the first one, and one can prove the second one similarly. Assume that (e) holds, then for any positive integer $m$ and any $w_0\in \mathfrak{h}^{0}$, we have
$$Z_R(\reg_{\shuffle}(S(y^m)\shuffle w_0-S(y^m)\ast w_0))=0.$$
Hence we get (v) by using Lemma \ref{Lem:Sym-ym}.
\qed

Finally, we recall two useful formulas to finish this section, which will be used later.

\begin{lem}[\cite{Ihara-Kajikawa-Ohno-Okuda,Ihara-Kaneko-Zagier}]
Let $k$ be a positive integer and $t$ be a variable. Then we have
\begin{align}
\exp_{\ast}\left(-\sum\limits_{n=1}^\infty\frac{z_{nk}}{n}t^n\right)=\frac{1}{1+z_kt}
\label{Eq:Exp-ast}
\end{align}
and
\begin{align}
\left(\frac{1}{1+z_kt}\right)\ast S\left(\frac{1}{1-z_kt}\right)=1.
\label{Eq:k-ast-Sk}
\end{align}
\end{lem}

Hence we have the following relations.

\begin{cor}
Assume that the $\mathbb{K}$-linear map $Z_R:\mathfrak{h}^0_{\ast}\rightarrow R$ is an algebra homomorphism. Then for any integer $k\geqslant 2$, we have
\begin{align}
\sum\limits_{n=0}^\infty\zeta_{R}(\{k\}^n)t^n=\exp\left(\sum\limits_{n=1}^\infty \frac{(-1)^{n-1}}{n}\zeta_{R}(nk)t^n\right)
\label{Eq:kn-nk}
\end{align}
and
\begin{align}
\left(\sum\limits_{n=0}^\infty(-1)^n\zeta_R(\{k\}^n)t^n\right)\left(\sum\limits_{n=0}^\infty\zeta_R^{\star}(\{k\}^n)t^n\right)=1,
\label{Eq:MZVk-MZSVk}
\end{align}
where $t$ is a variable.
\end{cor}

Here and below,  the notation $\{k_1,\ldots,k_r\}^n$ stands for $\underbrace{k_1,\ldots,k_r,\ldots,k_1,\ldots,k_r}_{nr \text{\;terms}}$.


\section{Relations deduced from regularized double shuffle relations}\label{Sec:Relations}

As stated in Section \ref{Sec:Intro}, it is conjectured that the regularized double shuffle relations give all algebraic relations among multiple zeta values, and hence one could obtain any algebraic relations from the regularized double shuffle relations. In this section, we reformulate and prove some of these results, which provides some evidence for the conjecture.

Below we always assume that the $\mathbb{K}$-linear map $Z_R:\mathfrak{h}^0\rightarrow R$ satisfies the regularized double shuffle relations unless otherwise specified.


\subsection{Linear double shuffle relations}

If the $\mathbb{K}$-linear map $Z_R:\mathfrak{h}^0\rightarrow R$ satisfies the double shuffle relations, we certainly obtain the following linear relations
$$Z_R(w_1\shuffle w_2-w_2\ast w_2)=0,\quad (\forall w_1,w_2\in\mathfrak{h}^0)$$
and
$$Z_R^{\star}(w_1\sshuffle w_2-w_2\sast w_2)=0,\quad (\forall w_1,w_2\in\mathfrak{h}^0),$$
which are equivalent and are called the linear double shuffle relations.


\subsection{Derivation relations}

As part of the regularized double shuffle relations, we have the following linear relations
$$Z_R(\partial_n(w_0))=0,\quad (\forall n\in\mathbb{N},\forall w_0\in\mathfrak{h}^0)$$
and
$$Z_R^{\star}(\partial_n^{\star}(w_0))=0,\quad (\forall n\in\mathbb{N},\forall w_0\in\mathfrak{h}^0),$$
which are equivalent and are called the derivation relations. The derivation relations for multiple zeta values were introduced by K. Ihara, M. Kaneko and D. Zagier in \cite{Ihara-Kaneko-Zagier}.


\subsection{Hoffman's relation}

In the algebra $\mathfrak{h}^1$, we have the following formulas.

\begin{lem}\label{Lem:Hoffman-Relation}
For any $w=z_{k_1}\cdots z_{k_n}\in\mathfrak{h}^1$ with $n,k_1,\ldots,k_n\in\mathbb{N}$, we have
\begin{align*}
\partial_1(w)=y\shuffle w-y\ast w=&\sum\limits_{i=1}^n\sum\limits_{j=0}^{k_i-2}z_{k_1}\cdots z_{k_{i-1}}z_{k_i-j}z_{j+1}z_{k_{i+1}}\cdots z_{k_n}\\
&-\sum\limits_{i=1}^nz_{k_1}\cdots z_{k_{i-1}}z_{k_i+1}z_{k_{i+1}}\cdots z_{k_n}
\end{align*}
and
\begin{align*}
\partial_1^{\star}(w)=y\sshuffle w-y\sast w=&\sum\limits_{i=1}^n\sum\limits_{j=0}^{k_i-2}z_{k_1}\cdots z_{k_{i-1}}z_{k_i-j}z_{j+1}z_{k_{i+1}}\cdots z_{k_n}\\
&-\sum\limits_{i=1}^n(k_i-1+\delta_{ni})z_{k_1}\cdots z_{k_{i-1}}z_{k_i+1}z_{k_{i+1}}\cdots z_{k_n},
\end{align*}
where $\delta_{ij}$ is Kronecker's delta symbol defined by
$$\delta_{ij}=\begin{cases}
1 & \text{if\;} i=j,\\
0 & \text{otherwise}.
\end{cases}$$
\end{lem}

\proof We have
\begin{align*}
\partial_1^{\star}(w)=&S^{-1}(\partial_1(S(w)))=S^{-1}(y\shuffle S(w)-y\ast S(w))\\
=&y\sshuffle w-y\sast w.
\end{align*}
Other equations can be obtained directly from the definitions.
\qed

From Lemma \ref{Lem:Hoffman-Relation}, we get the following Hoffman's relations
\begin{align*}
&\sum\limits_{i=1}^n\zeta_R(k_1,\ldots,k_{i-1},k_{i}+1,k_{i+1},\ldots,k_n)\\
=&\sum\limits_{i=1}^n\sum\limits_{j=0}^{k_i-2}\zeta_R(k_1,\ldots,k_{i-1},k_i-j,j+1,k_{i+1},\ldots,k_n)
\end{align*}
and
\begin{align*}
&\sum\limits_{i=1}^n(k_i-1+\delta_{ni})\zeta^{\star}_R(k_1,\ldots,k_{i-1},k_{i}+1,k_{i+1},\ldots,k_n)\\
=&\sum\limits_{i=1}^n\sum\limits_{j=0}^{k_i-2}\zeta^{\star}_R(k_1,\ldots,k_{i-1},k_i-j,j+1,k_{i+1},\ldots,k_n),
\end{align*}
where $n,k_1,\ldots,k_n\in\mathbb{N}$ with $k_1\geqslant 2$. The Hoffman's relation for multiple zeta values was first proved by M. E. Hoffman in \cite{Hoffman1992}. While the Hoffman's relation for multiple zeta-star values was first obtained by S. Muneta in \cite{Muneta} as a consequence of the regularized double shuffle relations.


\subsection{Sum formula}

For any integers $k,n$ with the conditions $k>n\geqslant 1$, we set
$$S(k,n)=x(x^{k-n-1}\shuffle y^{n-1})y.$$

\begin{lem}\label{Lem:Sum-Formula}
For any integers $k,n$ with the conditions $k>n\geqslant 1$, we have
\begin{align}
&\left(\overline{\sigma_{n-1}}-\sigma_{n-1}\right)(z_{k-n+1})=S(k,n)-z_k;
\label{Eq:Sigma}\\
&\left(\overline{\sigma_{n-1}}^{\star}-\sigma_{n-1}^{\star}\right)(z_{k-n+1})=\sum\limits_{i=1}^n(-1)^{n-i}\binom{k-i-1}{n-i}S(k,i)-z_k.
\label{Eq:Sigma-Star}
\end{align}
Moreover, if $n\geqslant 2$, we have
\begin{align}
(-1)^{n-1}\reg_{\shuffle}(y^{n-1}\ast z_{k-n+1})=S(k,n)-S(k,n-1).
\label{Eq:Shuffle-Reg}
\end{align}
\end{lem}

\proof Equation \eqref{Eq:Shuffle-Reg} is just \cite[Proposition 9]{Ihara-Kaneko-Zagier}. Since $\tau$ is an automorphism of $\mathfrak{h}_{\shuffle}$, we have
\begin{align*}
\overline{\sigma_{n-1}}(z_{k-n+1})=&\tau(\sigma_{n-1}(xy^{k-n}))=\tau(x(x^{n-1}\shuffle y^{k-n-1})y)\\
=&x(x^{k-n-1}\shuffle y^{n-1})y=S(k,n),
\end{align*}
which deduces equation \eqref{Eq:Sigma}. Now since $\sigma^{-1}$ is an automorphism of $\mathfrak{h}_{\shuffle}$, we get
\begin{align*}
\overline{\sigma_{n-1}}^{\star}(z_{k-n+1})=&S^{-1}(\overline{\sigma_{n-1}}(z_{k-n+1}))=S^{-1}(x(x^{k-n-1}\shuffle y^{n-1})y)\\
=&x(x^{k-n-1}\shuffle (-x+y)^{n-1})y\\
=&\sum\limits_{i=1}^{n}(-1)^{n-i}x(x^{k-n-1}\shuffle x^{n-i}\shuffle y^{i-1})y\\
=&\sum\limits_{i=1}^n(-1)^{n-i}\binom{k-i-1}{n-i}x(x^{k-i-1}\shuffle y^{i-1})y,
\end{align*}
which implies equation \eqref{Eq:Sigma-Star}.
\qed

Let $k$ and $n$ be integers with $k>n\geqslant 1$. From \eqref{Eq:Sigma-Star}, for any positive integer $m\leqslant n$, we have
$$\left(\overline{\sigma_{m-1}}^{\star}-\sigma_{m-1}^{\star}\right)(z_{k-m+1})=\sum\limits_{i=1}^m(-1)^{m-i}\binom{k-i-1}{m-i}S(k,i)-z_k.$$
Multiplying $\binom{k-m-1}{n-m}$ to the above equation and summing for $m$ from $1$ to $n$, we get
\begin{align}
\sum\limits_{m=1}^n\binom{k-m-1}{n-m}(\overline{\sigma_{m-1}}^{\star}-\sigma_{m-1}^{\star})(z_{k-m+1})=S(k,n)-\binom{k-1}{n-1}z_k.
\label{Eq:Sigma-Star-Sum}
\end{align}
Then by \eqref{Eq:Sigma} or \eqref{Eq:Shuffle-Reg}, we get the sum formula
\begin{align}
\sum\limits_{\wt(\mathbf{k})=k,\dep(\mathbf{k})=n\atop \mathbf{k}:\text{admissible}}\zeta_R(\mathbf{k})=\zeta_R(k);
\label{Eq:Sum-Formula}
\end{align}
and by \eqref{Eq:Sigma-Star-Sum}, we get the sum formula
\begin{align}
\sum\limits_{\wt(\mathbf{k})=k,\dep(\mathbf{k})=n\atop \mathbf{k}:\text{admissible}}\zeta^{\star}_R(\mathbf{k})=\binom{k-1}{n-1}\zeta_R(k).
\label{Eq:S-Sum-Formula}
\end{align}
The sum formula for multiple zeta values was first proved by A. Granville in \cite{Granville}, and it was shown by K. Ihara, M. Kaneko and D. Zagier in \cite{Ihara-Kaneko-Zagier} that one can deduce the sum formula for multiple zeta values from the regularized double shuffle relations. The equivalence of the sum formulas of multiple zeta values and of multiple zeta-star values was first proposed by M. E. Hoffman in \cite{Hoffman1992}.


\subsection{Weighted sum formula of L. Guo and B. Xie}

In this subsection, we discuss the weighted sum formula of L. Guo and B. Xie \cite{Guo-Xie2009}. They obtained their weighted sum formula just from the regularized double shuffle relations. Here we provide a much simpler proof of the weighted sum formula, and give the formula corresponding to the multiple zeta-star values.

For positive integers $n,k_1,\ldots,k_{n}$, let
\begin{align*}
\mathcal{C}(k_1,\ldots,k_{n})=&\sum\limits_{j=1}^{n}2^{k_1+\cdots+k_{j}-j}+2^{k_1+\cdots+k_{n}-n}\\
=&\sum\limits_{j=1}^{n-1}2^{k_1+\cdots+k_j-j}+2^{k_1+\cdots+k_{n}-n+1}\in\mathbb{N}.
\end{align*}

\begin{lem}[{\cite[Thoerem 2.4 and Theorem 2.6]{Guo-Xie2009}}]
For any integers $k,n$ with $k>n\geqslant 2$, we have
\begin{align}
&\sum\limits_{l+k_1+\cdots+k_{n-1}=k\atop l,k_i\geqslant 1, k_1\geqslant 2}z_l\ast z_{k_1}\cdots z_{k_{n-1}}=\sum\limits_{k_1+\cdots+k_n=k \atop k_i\geqslant 1,k_1=1,k_2\geqslant 2}z_{k_1}\cdots z_{k_n}\nonumber\\
&\qquad\qquad+(n-1)\sum\limits_{k_1+\cdots+k_n=k \atop k_i\geqslant 1,k_1\geqslant 2,k_2=1}z_{k_1}\cdots z_{k_n}+n\sum\limits_{k_1+\cdots+k_n=k \atop k_i\geqslant 1,k_1\geqslant 2,k_2\geqslant 2}z_{k_1}\cdots z_{k_n}\nonumber\\
&\qquad\qquad+(k-n)\sum\limits_{k_1+\cdots+k_{n-1}=k \atop k_i\geqslant 1,k_1\geqslant 2}z_{k_1}\cdots z_{k_{n-1}}
\label{Eq:Ast-1-n-Admissible}
\end{align}
and
\begin{align}
&\sum\limits_{l+k_1+\cdots+k_{n-1}=k\atop l,k_i\geqslant 1,k_1\geqslant 2}z_l\shuffle z_{k_1}\cdots z_{k_{n-1}}
=\sum\limits_{k_1+\cdots+k_n=k\atop k_i\geqslant 1}[\mathcal{C}(k_1,\ldots,k_{n-1})\nonumber\\
&\qquad\qquad-\mathcal{C}(k_2,\ldots,k_{n-1})]z_{k_1}\cdots z_{k_n}-\sum\limits_{k_1+\cdots+k_{n}=k\atop k_i\geqslant 1, k_2=1}z_{k_1}\cdots z_{k_{n}}.
\label{Eq:Shuffle-1-n-Admissible}
\end{align}
Here, if $n=2$, we set $\mathcal{C}(k_2,\ldots,k_{n-1})=1$.
\end{lem}

\proof In \cite{Guo-Xie2009}, L. Guo and B. Xie proved these formulas by induction on $n$. Here we provide a much simpler proof. For \eqref{Eq:Ast-1-n-Admissible}, we use the formula
\begin{align*}
z_l\ast z_{k_1}\cdots z_{k_{n-1}}=&\sum\limits_{i=0}^{n-1}z_{k_1}\cdots z_{k_i}z_lz_{k_{i+1}}\cdots z_{k_{n-1}}\\
&+\sum\limits_{i=1}^{n-1}z_{k_1}\cdots z_{k_{i-1}}z_{l+k_i}z_{k_{i+1}}\cdots z_{k_{n-1}},
\end{align*}
which is immediate from the combinatorial description of the product $\ast$. Hence the left-hand side of \eqref{Eq:Ast-1-n-Admissible} is
\begin{align*}
&\sum\limits_{l+k_1+\cdots+k_{n-1}=k\atop l,k_j\geqslant 1,k_1\geqslant 2}\sum\limits_{i=0}^{n-1}z_{k_1}\cdots z_{k_i}z_lz_{k_{i+1}}\cdots z_{k_{n-1}}\\
&+\sum\limits_{l+k_1+\cdots+k_{n-1}=k\atop l,k_j\geqslant 1,k_1\geqslant 2}\sum\limits_{i=1}^{n-1}z_{k_1}\cdots z_{k_{i-1}}z_{l+k_i}z_{k_{i+1}}\cdots z_{k_{n-1}}\\
=&\sum\limits_{k_1+\cdots+k_n=k\atop k_i\geqslant 1,k_2\geqslant 2}z_{k_1}\cdots z_{k_n}+(n-1)\sum\limits_{k_1+\cdots+k_n=k\atop k_i\geqslant 1,k_1\geqslant 2}z_{k_1}\cdots z_{k_n}\\
&+\sum\limits_{k_1+\cdots+k_{n-1}=k\atop k_i\geqslant 1,k_1\geqslant 2}(k_1-2)z_{k_1}\cdots z_{k_{n-1}}+\sum\limits_{i=2}^{n-1}\sum\limits_{k_1+\cdots+k_{n-1}=k\atop k_j\geqslant 1,k_1\geqslant 2}(k_i-1)z_{k_1}\cdots z_{k_{n-1}},
\end{align*}
which is just the right-hand side of \eqref{Eq:Ast-1-n-Admissible}.

For \eqref{Eq:Shuffle-1-n-Admissible}, we use the formula
\begin{align*}
&z_l\shuffle z_{k_1}\cdots z_{k_{n-1}}\\
=&\sum\limits_{i=1}^{n-1}\sum\limits_{\alpha_j\geqslant 1,\alpha_1+\cdots+\alpha_{i+1}\atop =l+k_1+\cdots+k_{i}}\prod\limits_{j=1}^{i-1}\binom{\alpha_j-1}{k_j-1}
\binom{\alpha_i-1}{k_i-\alpha_{i+1}}
z_{\alpha_1}\cdots z_{\alpha_{i+1}}z_{k_{i+1}}\cdots z_{k_{n-1}}\\
&+\sum\limits_{\alpha_j\geqslant 1,\alpha_1+\cdots+\alpha_n\atop =l+k_1+\cdots+k_{n-1}}\prod\limits_{j=1}^{n-1}\binom{\alpha_j-1}{k_j-1}z_{\alpha_1}\cdots z_{\alpha_n},
\end{align*}
which is \cite[Proposition 2.7]{Li-Qin2017} and can be proved by the combinatorial description of the product $\shuffle$. Hence the left-hand side of \eqref{Eq:Shuffle-1-n-Admissible} is
\begin{align*}
&\sum\limits_{\alpha_1+\alpha_2+k_2+\cdots+k_{n-1}=k\atop \alpha_i\geqslant 1,k_j\geqslant 1}\sum\limits_{k_1\geqslant 2,k_1\geqslant \alpha_2}\binom{\alpha_1-1}{k_1-\alpha_2}z_{\alpha_1}z_{\alpha_2}z_{k_2}\cdots z_{k_{n-1}}\\
&+\sum\limits_{i=2}^{n-1}\sum\limits_{\alpha_1+\cdots+\alpha_{i+1}+k_{i+1}+\cdots+k_{n-1}=k\atop \alpha_s\geqslant 1,k_t\geqslant 1}\sum\limits_{k_1=2}^{\alpha_1}\binom{\alpha_1-1}{k_1-1}\prod\limits_{j=2}^{i-1}\sum\limits_{k_j=1}^{\alpha_j}\binom{\alpha_j-1}{k_j-1}
\\
&\times \sum\limits_{k_i=\alpha_{i+1}}^{\alpha_i+\alpha_{i+1}-1}\binom{\alpha_i-1}{k_i-\alpha_{i+1}}z_{\alpha_1}\cdots z_{\alpha_{i+1}}z_{k_{i+1}}\cdots z_{k_{n-1}}\\
&+\sum\limits_{\alpha_1+\cdots+\alpha_n=k\atop \alpha_i\geqslant 1}\sum\limits_{k_1=2}^{\alpha_1}\binom{\alpha_1-1}{k_1-1}\prod\limits_{j=2}^{n-1}\sum\limits_{k_j=1}^{\alpha_j}\binom{\alpha_j-1}{k_j-1}z_{\alpha_1}\cdots z_{\alpha_n}\\
=&\sum\limits_{k_1+\cdots+k_n=k\atop k_i\geqslant 1,k_2\geqslant 2}2^{k_1-1}z_{k_1}\cdots z_{k_n}+\sum\limits_{k_1+\cdots+k_n=k\atop k_i\geqslant 1,k_2=1}\left(2^{k_1-1}-1\right)z_{k_1}\cdots z_{k_n}\\
&+\sum\limits_{i=2}^{n-1}\sum\limits_{k_1+\cdots+k_n=k\atop k_j\geqslant 1}\left(2^{k_1+\cdots+k_i-i}-2^{k_2+\cdots+k_{i}-(i-1)}\right)z_{k_1}\cdots z_{k_n}\\
&+\sum\limits_{k_1+\cdots+k_n=k\atop k_i\geqslant 1}\left(2^{k_1+\cdots+k_{n-1}-(n-1)}-2^{k_2+\cdots+k_{n-1}-(n-2)}\right)z_{k_1}\cdots z_{k_n},
\end{align*}
which is just the right-hand side of \eqref{Eq:Shuffle-1-n-Admissible}.
\qed

We can obtain the similar formulas for the products $\sast$ and $\sshuffle$.

\begin{lem}
For any integers $k,n$ with $k>n\geqslant 2$, we have
\begin{align}
&\sum\limits_{l+k_1+\cdots+k_{n-1}=k\atop l,k_i\geqslant 1,k_1\geqslant 2}z_l\sast z_{k_1}\cdots z_{k_{n-1}}=\sum\limits_{k_1+\cdots+k_n=k\atop k_i\geqslant 1,k_1=1,k_2\geqslant 2}z_{k_1}\cdots z_{k_n}\nonumber\\
&\qquad\qquad+(n-1)\sum\limits_{k_1+\cdots+k_n=k\atop k_i\geqslant 1,k_1\geqslant 2,k_2=1}z_{k_1}\cdots z_{k_n}+n\sum\limits_{k_1+\cdots+k_n=k\atop k_i\geqslant 1,k_1\geqslant 2,k_2\geqslant 2}z_{k_1}\cdots z_{k_n}\nonumber\\
&\qquad\qquad-(k-n)\sum\limits_{k_1+\cdots+k_{n-1}=k\atop k_i\geqslant 1,k_1\geqslant 2}z_{k_1}\cdots z_{k_{n-1}}
\label{Eq:S-Ast-1-n-Admissible}
\end{align}
and
\begin{align}
&\sum\limits_{l+k_1+\cdots+k_{n-1}=k\atop l,k_i\geqslant 1,k_1\geqslant 2}z_l\sshuffle z_{k_1}\cdots z_{k_{n-1}}\nonumber\\
=&-\sum\limits_{k_1+\cdots+k_n=k\atop k_i\geqslant 1, k_2=1}z_{k_1}\cdots z_{k_n}+\sum\limits_{k_1+\cdots+k_{n-1}=k\atop k_i\geqslant 1,k_1\geqslant 2}z_{k_1}\cdots z_{k_{n-1}}\nonumber\\
&+\sum\limits_{k_1+\cdots+k_n=k\atop k_i\geqslant 1}\left[\mathcal{C}(k_1,\ldots,k_{n-1})-\mathcal{C}(k_2,\ldots,k_{n-1})\right]z_{k_1}\cdots z_{k_n}\nonumber\\
&-\sum\limits_{k_1+\cdots+k_{n-1}=k\atop k_i\geqslant 1}\left[\mathcal{C}(k_1,\ldots,k_{n-1})-\mathcal{C}(k_2,\ldots,k_{n-1})\right.\nonumber\\
&\left.-\mathcal{C}(k_1,\ldots,k_{n-2})
+\mathcal{C}(k_2,\ldots,k_{n-2})\right]z_{k_1}\cdots z_{k_{n-1}}.
\label{Eq:S-Shuffle-1-n-Admissible}
\end{align}
Here, for $n=2$, we set $\mathcal{C}(k_1,\ldots,k_{n-2})=\mathcal{C}(k_2,\ldots,k_{n-1})=1$, $\mathcal{C}(k_2,\ldots,k_{n-2})=1-k$; and for $n=3$, we set $\mathcal{C}(k_2,\ldots,k_{n-2})=1$.
\end{lem}

\proof Similar as \eqref{Eq:Ast-1-n-Admissible}, we use the formula
\begin{align*}
z_l\sast z_{k_1}\cdots z_{k_{n-1}}=&\sum\limits_{i=0}^{n-1}z_{k_1}\cdots z_{k_i}z_lz_{k_{i+1}}\cdots z_{k_{n-1}}\\
&-\sum\limits_{i=1}^{n-1}z_{k_1}\cdots z_{k_{i-1}}z_{l+k_i}z_{k_{i+1}}\cdots z_{k_{n-1}}
\end{align*}
to prove \eqref{Eq:S-Ast-1-n-Admissible}. For \eqref{Eq:S-Shuffle-1-n-Admissible}, we use the formula
\begin{align*}
&z_l\sshuffle z_{k_1}\cdots z_{k_{n-1}}=\sum\limits_{i=1}^{n-1}\sum\limits_{\alpha_j\geqslant 1,\alpha_1+\cdots+\alpha_{i+1}\atop =l+k_1+\cdots+k_{i}}\prod\limits_{j=1}^{i-1}\binom{\alpha_j-1}{k_j-1}\binom{\alpha_i-1}{k_i-\alpha_{i+1}}\\
&\times (z_{\alpha_1}\cdots z_{\alpha_{i+1}}z_{k_{i+1}}\cdots z_{k_{n-1}}-z_{\alpha_1}\cdots z_{\alpha_{i-1}}z_{\alpha_i+\alpha_{i+1}}z_{k_{i+1}}\cdots z_{k_{n-1}})\\
&+\sum\limits_{\alpha_j\geqslant 1,\alpha_1+\cdots+\alpha_n\atop =l+k_1+\cdots+k_{n-1}}\prod\limits_{j=1}^{n-1}\binom{\alpha_j-1}{k_j-1}(z_{\alpha_1}\cdots z_{\alpha_n}-z_{\alpha_1}\cdots z_{\alpha_{n-2}}z_{\alpha_{n-1}+\alpha_n}).
\end{align*}
Hence similar as the proof of \eqref{Eq:Shuffle-1-n-Admissible}, the left-hand side of \eqref{Eq:S-Shuffle-1-n-Admissible} is
\begin{align*}
&\sum\limits_{k_1+\cdots+k_n=k\atop k_i\geqslant 1}\left[\mathcal{C}(k_1,\ldots,k_{n-1})-\mathcal{C}(k_2,\ldots,k_{n-1})\right]z_{k_1}\cdots z_{k_n}-\sum\limits_{k_1+\cdots+k_{n}=k\atop k_i\geqslant 1,k_2=1}z_{k_1}\cdots z_{k_{n}}\\
&-\sum\limits_{\alpha_1+\alpha_2+k_2+\cdots+k_{n-1}=k\atop \alpha_i,k_j\geqslant 1}\sum\limits_{k_1\geqslant 2,k_1\geqslant \alpha_2}\binom{\alpha_1-1}{k_1-\alpha_2}z_{\alpha_1+\alpha_2}z_{k_2}\cdots z_{k_{n-1}}\\
&-\sum\limits_{i=2}^{n-1}\sum\limits_{k_1+\cdots+k_{n}=k\atop k_j\geqslant 1}\left[2^{k_1+\cdots+k_i-i}-2^{k_2+\cdots+k_i-(i-1)}\right]z_{k_1}\cdots z_{k_{i-1}}z_{k_i+k_{i+1}}z_{k_{i+2}}\cdots z_{k_n}\\
&-\sum\limits_{k_1+\cdots+k_n=k\atop k_i\geqslant 1}\left[2^{k_1+\cdots+k_{n-1}-(n-1)}-2^{k_2+\cdots+k_{n-1}-(n-2)}\right]z_{k_1}\cdots z_{k_{n-2}}z_{k_{n-1}+k_n}.
\end{align*}
Denote the last three terms by $\Sigma_1,\Sigma_2$ and $\Sigma_3$, respectively. We find $-\Sigma_1$ is
\begin{align*}
&\sum\limits_{k_1+\cdots+k_n=k\atop k_i\geqslant 1,k_2\geqslant 2}2^{k_1-1}z_{k_1+k_2}z_{k_3}\cdots z_{k_n}+\sum\limits_{k_1+\cdots+k_n=k\atop k_i\geqslant 1,k_2=1}(2^{k_1-1}-1)z_{k_1+k_2}z_{k_3}\cdots z_{k_n}\\
=&\sum\limits_{k_1+\cdots+k_n=k\atop k_i\geqslant 1}2^{k_1-1}z_{k_1+k_2}z_{k_3}\cdots z_{k_n}-\sum\limits_{k_1+\cdots+k_n=k\atop k_i\geqslant 1,k_2=1}z_{k_1+k_2}z_{k_3}\cdots z_{k_n}\\
=&\sum\limits_{k_1+\cdots+k_{n-1}=k\atop k_i\geqslant 1}(2^{k_1-1}-1)z_{k_1}\cdots z_{k_{n-1}}-\sum\limits_{k_1+\cdots+k_{n-1}=k\atop k_i\geqslant 1,k_1\geqslant 2}z_{k_1}\cdots z_{k_{n-1}}.
\end{align*}
Similarly, $-\Sigma_2$ is
\begin{align*}
&\sum\limits_{i=2}^{n-1}\sum\limits_{k_1+\cdots+k_{n-1}=k\atop k_j\geqslant 1}\left[2^{k_1+\cdots+k_{i-1}-i}-2^{k_2+\cdots+k_{i-1}-(i-1)}\right](2^{k_i}-2)z_{k_1}\cdots z_{k_{n-1}}\\
=&\sum\limits_{i=2}^{n-1}\sum\limits_{k_1+\cdots+k_{n-1}=k\atop k_j\geqslant 1}\left[2^{k_1+\cdots+k_i-i}-2^{k_2+\cdots+k_i-(i-1)}-2^{k_1+\cdots+k_{i-1}-(i-1)}
\right.\\
&\left.+2^{k_2+\cdots+k_{i-1}-(i-2)}\right]z_{k_1}\cdots z_{k_{n-1}}.
\end{align*}
And finally, $-\Sigma_3$ is
\begin{align*}
&\sum\limits_{k_1+\cdots+k_{n-1}=k\atop k_i\geqslant 1}\left[2^{k_1+\cdots+k_{n-1}-(n-1)}-2^{k_2+\cdots+k_{n-1}-(n-2)}-2^{k_1+\cdots+k_{n-2}-(n-2)}
\right.\\
&\left.+2^{k_2+\cdots+k_{n-2}-(n-3)}\right]z_{k_1}\cdots z_{k_{n-1}}.
\end{align*}
Thus $-\Sigma_1-\Sigma_2-\Sigma_3$ is
\begin{align*}
&\sum\limits_{k_1+\cdots+k_{n-1}=k\atop k_i\geqslant 1}\left[\mathcal{C}(k_1,\ldots,k_{n-1})-\mathcal{C}(k_2,\ldots,k_{n-1})-\mathcal{C}(k_1,\ldots,k_{n-2})\right.\\
&\left.+\mathcal{C}(k_2,\ldots,k_{n-2})\right]z_{k_1}\cdots z_{k_{n-1}}-\sum\limits_{k_1+\cdots+k_{n-1}=k\atop k_i\geqslant 1,k_1\geqslant 2}z_{k_1}\cdots z_{k_{n-1}}.
\end{align*}
Now \eqref{Eq:S-Shuffle-1-n-Admissible} follows immediately.
\qed

Note that if $n\geqslant 3$, for any positive integers $k_1,\ldots,k_{n-1}$, we have
\begin{align*}
&\mathcal{C}(k_1,\ldots,k_{n-1})-\mathcal{C}(k_2,\ldots,k_{n-1})-\mathcal{C}(k_1,\ldots,k_{n-2})
+\mathcal{C}(k_2,\ldots,k_{n-2})\\
=&2^{\sum\limits_{j=2}^{n-2}k_j-n+3}(2^{k_1-1}-1)(2^{k_{n-1}}-1).
\end{align*}

From \eqref{Eq:Ast-1-n-Admissible}-\eqref{Eq:S-Shuffle-1-n-Admissible}, we immediately get the following result.

\begin{cor}
For any integers $k,n$ with $k>n\geqslant 2$, we have
\begin{align}
&\sum\limits_{l+k_1+\cdots+k_{n-1}=k\atop l,k_i\geqslant 1,k_1\geqslant 2}(z_l\shuffle z_{k_1}\cdots z_{k_{n-1}}-z_l\ast z_{k_1}\cdots z_{k_{n-1}})\nonumber\\
=&\sum\limits_{k_1+\cdots+k_n=k\atop k_i\geqslant 1,k_1\geqslant 2}[\mathcal{C}(k_1,\ldots,k_{n-1})-\mathcal{C}(k_2,\ldots,k_{n-1})]z_{k_1}\cdots z_{k_n}\nonumber\\
&-n\sum\limits_{k_1+\cdots+k_n=k\atop k_i\geqslant 1,k_1\geqslant 2}z_{k_1}\cdots z_{k_n}-(k-n)\sum\limits_{k_1+\cdots+k_{n-1}=k\atop k_i\geqslant 1,k_1\geqslant 2}z_{k_1}\cdots z_{k_{n-1}}
\label{Eq:WSum-GuoXie}
\end{align}
and
\begin{align}
&\sum\limits_{l+k_1+\cdots+k_{n-1}=k\atop l,k_i\geqslant 1,k_1\geqslant 2}(z_l\sshuffle z_{k_1}\cdots z_{k_{n-1}}-z_l\sast z_{k_1}\cdots z_{k_{n-1}})\nonumber\\
=&\sum\limits_{k_1+\cdots+k_n=k\atop k_i\geqslant 1,k_1\geqslant 2}\left[\mathcal{C}(k_1,\ldots,k_{n-1})-\mathcal{C}(k_2,\ldots,k_{n-1})\right]z_{k_1}\cdots z_{k_n}\nonumber\\
&-\sum\limits_{k_1+\cdots+k_{n-1}=k\atop k_i\geqslant 1,k_1\geqslant 2}\left[\mathcal{C}(k_1,\ldots,k_{n-1})-\mathcal{C}(k_2,\ldots,k_{n-1})\right.\nonumber\\
&\left.-\mathcal{C}(k_1,\ldots,k_{n-2})
+\mathcal{C}(k_2,\ldots,k_{n-2})\right]z_{k_1}\cdots z_{k_{n-1}}\nonumber\\
&-n\sum\limits_{k_1+\cdots+k_n=k\atop k_i\geqslant 1,k_1\geqslant 2}z_{k_1}\cdots z_{k_n}+(k-n+1)\sum\limits_{k_1+\cdots+k_{n-1}=k\atop k_i\geqslant 1,k_1\geqslant 2}z_{k_1}\cdots z_{k_{n-1}}.
\label{Eq:S-WSum-GuoXie}
\end{align}
\end{cor}

Applying the map $Z_R$ to \eqref{Eq:WSum-GuoXie} and the map $Z_R^{\star}$ to \eqref{Eq:S-WSum-GuoXie}, and using the sum formulas \eqref{Eq:Sum-Formula} and \eqref{Eq:S-Sum-Formula}, we get the weighted sum formulas
\begin{align}
\sum\limits_{k_1+\cdots+k_n=k\atop k_i\geqslant 1,k_1\geqslant 2}[\mathcal{C}(k_1,\ldots,k_{n-1})-\mathcal{C}(k_2,\ldots,k_{n-2})]\zeta_R(k_1,\ldots,k_n)=k\zeta_R(k)
\label{Eq:WSum-GuoXie-MZV}
\end{align}
and
\begin{align*}
&\sum\limits_{k_1+\cdots+k_n=k\atop k_i\geqslant 1,k_1\geqslant 2}\left[\mathcal{C}(k_1,\ldots,k_{n-1})-\mathcal{C}(k_2,\ldots,k_{n-1})\right]\zeta_R^{\star}(k_1,\ldots,k_n)\\
&-\sum\limits_{k_1+\cdots+k_{n-1}=k\atop k_i\geqslant 1,k_1\geqslant 2}\left[\mathcal{C}(k_1,\ldots,k_{n-1})-\mathcal{C}(k_2,\ldots,k_{n-1})\right.\\
&\left.-\mathcal{C}(k_1,\ldots,k_{n-2})+\mathcal{C}(k_2,\ldots,k_{n-2})\right]\zeta_R^{\star}(k_1,\ldots,k_{n-1})=\binom{k-1}{n-1}\zeta_R(k),
\end{align*}
where $k$ and $n$ are integers with $k>n\geqslant 2$. Note that if $Z_R=Z$, then \eqref{Eq:WSum-GuoXie-MZV} is \cite[Theorem 2.2]{Guo-Xie2009}. Taking $n=2$ and using the sum formulas, we find for any $k\geqslant 3$,
\begin{align*}
&\sum\limits_{i=2}^{k-1}2^i\zeta_R(i,k-i)=(k+1)\zeta_R(k),\\
&\sum\limits_{i=2}^{k-1}2^i\zeta_R^{\star}(i,k-i)=(2^k+k-3)\zeta_R(k),
\end{align*}
which were first discovered by Y. Ohno and W. Zudilin in \cite{Ohno-Zudilin}.


\subsection{Restricted sum formula of M. Eie, W. C. Liaw and Y. L. Ong}

We compute explicitly the relations coming from
$$Z_R((\sigma_m-\overline{\sigma_m})(w_0))=0$$
for any nonempty word $w_0\in\mathfrak{h}^0$ and any nonnegative integer $m$.

\begin{lem}
Let $w_0=x^{a_1}y^{b_1}\cdots x^{a_h}y^{b_h}\in\mathfrak{h}^0$ with $h$, $a_1,\ldots,a_h$, $b_1,\ldots,b_h\in\mathbb{N}$. Then for any nonnegative integer $m$, we have
\begin{align}
\sigma_m(w_0)=\sum\limits_{\varepsilon_1+\cdots+\varepsilon_h=m\atop \varepsilon_i\geqslant 0}x^{a_1}(x^{\varepsilon_1}\shuffle y^{b_1-1})yx^{a_2}(x^{\varepsilon_2}\shuffle y^{b_2-1})y\cdots x^{a_h}(x^{\varepsilon_h}\shuffle y^{b_h-1})y
\label{Eq:Sigma-m}
\end{align}
and
\begin{align}
\overline{\sigma_m}(w_0)=\sum\limits_{\varepsilon_1+\cdots+\varepsilon_h=m\atop \varepsilon_i\geqslant 0}x(x^{a_1-1}\shuffle y^{\varepsilon_1})y^{b_1}x(x^{a_2-1}\shuffle y^{\varepsilon_2})y^{b_2}\cdots x(x^{a_h-1}\shuffle y^{\varepsilon_h})y^{b_h}.
\label{Eq:Sigma-bar-m}
\end{align}
\end{lem}

\proof Since $w_0=z_{a_1+1}z_1^{b_1-1}z_{a_2+1}z_1^{b_2-1}\cdots z_{a_h+1}z_1^{b_h-1}$, we have
\begin{align*}
\sigma_m(w_0)=&\sum\limits_{\varepsilon_1+\cdots+\varepsilon_h=m\atop \varepsilon_i\geqslant 0}\sum\limits_{\varepsilon_{i1}+\cdots+\varepsilon_{ib_i}=\varepsilon_i\atop \varepsilon_{ij}\geqslant 0}\prod\limits_{i=1}^hz_{a_i+\varepsilon_{i1}+1}z_{\varepsilon_{i2}+1}\cdots z_{\varepsilon_{ib_i}+1}\\
=&\sum\limits_{\varepsilon_1+\cdots+\varepsilon_h=m\atop \varepsilon_i\geqslant 0}\prod\limits_{i=1}^hx^{a_i}\left(\sum\limits_{\varepsilon_{i1}+\cdots+\varepsilon_{ib_i}=\varepsilon_i\atop \varepsilon_{ij}\geqslant 0}x^{\varepsilon_{i1}}yx^{\varepsilon_{i2}}y\cdots x^{\varepsilon_{ib_i}}\right)y\\
=&\sum\limits_{\varepsilon_1+\cdots+\varepsilon_h=m\atop \varepsilon_i\geqslant 0}\prod\limits_{i=1}^hx^{a_i}(x^{\varepsilon_i}\shuffle y^{b_i-1})y,
\end{align*}
which is \eqref{Eq:Sigma-m}. The equation \eqref{Eq:Sigma-bar-m} follows from \eqref{Eq:Sigma-m} and the fact that $\tau$ is an automorphism of the shuffle algebra $\mathfrak{h}_{\shuffle}$.
\qed

With the help of \eqref{Eq:Sigma-m} and \eqref{Eq:Sigma-bar-m},  for any positive integers $h$, $a_1$, $\ldots$, $a_h$, $b_1$, $\ldots$, $b_h$ and any nonnegative integers $m$, we obtain the relation
\begin{align}
&\sum\limits_{p_1+\cdots+p_h=m\atop p_i\geqslant 0}\sum\limits_{k_{i1}+\cdots+k_{i,p_i+1}=a_i+p_i \atop k_{ij}\geqslant 1}\zeta_R(k_{11}+1,k_{12},\ldots,k_{1,p_1+1},\{1\}^{b_1-1},\nonumber\\
&\qquad\qquad\qquad\qquad\qquad\ldots,k_{h1}+1,k_{h2},\ldots,k_{h,p_h+1},\{1\}^{b_h-1})\nonumber\\
=&\sum\limits_{q_1+\cdots+q_h=m\atop q_i\geqslant 0}\sum\limits_{l_{i1}+\cdots+l_{ib_i}=b_i+q_i\atop l_{ij}\geqslant 1}\zeta_R(a_1+l_{11},l_{12},\ldots,l_{1b_1},\ldots,a_h+l_{h1},l_{h2},\ldots,l_{hb_h}),
\label{Eq:Ohno-OneIndex}
\end{align}
which is called the vector version of the restricted sum formula in \cite[Theorem 3.4.1]{Eie-book}. If $s=1$, we get the restricted sum formula (\cite{Eie-Liaw-Ong}): for any positive integers $a,b$ and any nonnegative integer $m$, we have
$$\sum\limits_{k_1+\cdots+k_{m+1}=a+m\atop k_i\geqslant 1}\zeta_R(k_1+1,k_2,\ldots,k_{m+1},\{1\}^{b-1})=\sum\limits_{l_1+\cdots+l_b=b+m\atop l_i\geqslant 1}\zeta_R(a+l_1,l_2,\ldots,l_b),$$
which is the sum formula \eqref{Eq:Sum-Formula} in the case of $b=1$.

Now we consider the relations coming from
\begin{align}
Z_R^{\star}((\sigma_m^{\star}-\overline{\sigma_m}^{\star})(w_0))=0
\label{Eq:sigma-tau-sigma-S}
\end{align}
for any nonempty word $w_0\in\mathfrak{h}^0$  and any nonnegative integer $m$. Since $S^{-1}:\mathfrak{h}^0\rightarrow\mathfrak{h}^0$ is  invertible, we can  take $w_0=S^{-1}(x^{a_1}y^{b_1}\cdots x^{a_h}y^{b_h})$, where $h,a_1,b_1,\ldots,a_h,b_h$ are any positive integers.

\begin{lem}
Let $w_0=S^{-1}(x^{a_1}y^{b_1}\cdots x^{a_h}y^{b_h})$ with $h,a_1,b_1,\ldots,a_h,b_h\in\mathbb{N}$. Then we have
\begin{align}
&\sigma_m^{\star}(w_0)=\sum\limits_{\varepsilon_1+\cdots+\varepsilon_h=m\atop \varepsilon_i\geqslant 0}\sum\limits_{0\leqslant j_i\leqslant b_i-1\atop i=1,\ldots,h}(-1)^{\sum\limits_{l=1}^h(b_l-j_l-1)}\prod\limits_{l=1}^h\binom{\varepsilon_l+b_l-j_l-1}{\varepsilon_l}\nonumber\\
&\quad \times \prod\limits_{l=1}^{h-1}[x^{a_l}(x^{\varepsilon_l+b_l-j_l-1}\shuffle y^{j_l})(-x+y)]x^{a_h}(x^{\varepsilon_h+b_h-j_h-1}\shuffle y^{j_h})y
\label{Eq:S-Sigma-m}
\end{align}
and
\begin{align}
&\overline{\sigma_m}^{\star}(w_0)=\sum\limits_{\varepsilon_1+\cdots+\varepsilon_h=m\atop \varepsilon_i\geqslant 0}\sum\limits_{0\leqslant j_i\leqslant \varepsilon_i\atop i=1,\ldots,h}(-1)^{\sum\limits_{l=1}^h(\varepsilon_l-j_l)}\prod\limits_{l=1}^h\binom{a_l+\varepsilon_l-j_l-1}{a_l-1}\nonumber\\
&\times \prod\limits_{l=1}^{h-1}[x(x^{a_l+\varepsilon_l-j_l-1}\shuffle y^{j_l})(-x+y)^{b_l}]x(x^{a_h+\varepsilon_h-j_h-1}\shuffle y^{j_h})(-x+y)^{b_h-1}y.
\label{Eq:S-Sigma-bar-m}
\end{align}
\end{lem}

\proof By \eqref{Eq:Sigma-m}, we get
$$\sigma_m^{\star}(w_0)=\sum\limits_{\varepsilon_1+\cdots+\varepsilon_h=m\atop \varepsilon_i\geqslant 0}S^{-1}\left(x^{a_1}(x^{\varepsilon_1}\shuffle y^{b_1-1})y\cdots x^{a_h}(x^{\varepsilon_h}\shuffle y^{b_h-1})y\right).$$
Since $\sigma^{-1}$ is an automorphism of $\mathfrak{h}_{\shuffle}$, we find $\sigma_m^{\star}(w_0)$ is
$$\sum\limits_{\varepsilon_1+\cdots+\varepsilon_h=m\atop \varepsilon_i\geqslant 0}\prod\limits_{i=1}^{h-1}\left[x^{a_i}(x^{\varepsilon_i}\shuffle (-x+y)^{b_i-1})(-x+y)\right]x^{a_h}(x^{\varepsilon_h}\shuffle (-x+y)^{b_h-1})y,$$
from which we get \eqref{Eq:S-Sigma-m}. Similarly, one can prove \eqref{Eq:S-Sigma-bar-m}.
\qed

From \eqref{Eq:S-Sigma-m} and \eqref{Eq:S-Sigma-bar-m},  for any positive integers $h$, $a_1,\ldots,a_h$, $b_1,\ldots,b_h$ and any nonnegative integers $m$, we have the relation
\begin{align*}
&\sum\limits_{\varepsilon_1+\cdots+\varepsilon_h=m\atop \varepsilon_i\geqslant 0}\sum\limits_{0\leqslant j_i\leqslant b_i-1\atop i=1,\ldots,h}(-1)^{\sum\limits_{l=1}^h(b_l-j_l-1)}\prod\limits_{l=1}^h\binom{\varepsilon_l+b_l-j_l-1}{\varepsilon_l}\\
&\quad \times Z_R^{\star}\left(\prod\limits_{l=1}^{h-1}[x^{a_l}(x^{\varepsilon_l+b_l-j_l-1}\shuffle y^{j_l})(-x+y)]x^{a_h}(x^{\varepsilon_h+b_h-j_h-1}\shuffle y^{j_h})y\right)\\
=&\sum\limits_{\varepsilon_1+\cdots+\varepsilon_h=m\atop \varepsilon_i\geqslant 0}\sum\limits_{0\leqslant j_i\leqslant \varepsilon_i\atop i=1,\ldots,h}(-1)^{\sum\limits_{l=1}^h(\varepsilon_l-j_l)}\prod\limits_{l=1}^h\binom{a_l+\varepsilon_l-j_l-1}{a_l-1}\\
&\times Z_R^{\star}\left(\prod\limits_{l=1}^{h-1}[x(x^{a_l+\varepsilon_l-j_l-1}\shuffle y^{j_l})(-x+y)^{b_l}]x(x^{a_h+\varepsilon_h-j_h-1}\shuffle y^{j_h})(-x+y)^{b_h-1}y\right).
\end{align*}
Let $h=1$. Then for any positive integers $a,b$ and any nonnegative integer $m$, we have
\begin{align*}
&\sum\limits_{j=0}^{b-1}(-1)^{j}\binom{m+b-j-1}{m}\sum\limits_{k_1+\cdots+k_{j+1}=m+b\atop k_i\geqslant 1}\zeta_R^{\star}(a+k_1,k_2,\ldots,k_{j+1}+\delta_{j0}a)\\
=&\sum\limits_{j=0}^m\sum\limits_{l=0}^{b-1}(-1)^{m-j-l}\binom{a+m-j-1}{a-1}\sum\limits_{{k_i\geqslant 1,k_1+\cdots+k_{j+1}=a+m\atop n_i\geqslant 1,n_1+\cdots+n_{l+1}=b}}\\
&\quad\times \zeta_R^{\star}(k_1+1,k_2,\ldots,k_j,k_{j+1}+n_1-1+\delta_{j0},n_2,\ldots,n_{l+1}).
\end{align*}
Moreover, let $b=1$ in the above formula. Then for any integers $k,n$ with $k>n\geqslant 1$, we have
\begin{align*}
\sum\limits_{j=1}^n(-1)^{n-j}\binom{k-j-1}{k-n-1}\sum\limits_{\wt(\mathbf{k})=k,\dep(\mathbf{k})=j\atop \mathbf{k}:\text{admissible}}\zeta_R^{\star}(\mathbf{k})=\zeta_R(k),
\end{align*}
which is in fact equivalent to the sum formula \eqref{Eq:S-Sum-Formula}.


\subsection{Restricted sum formulas of double zeta values}

In \cite{Gangl-Kaneko-Zagier}, it was proved by H. Gangl, M. Kaneko and D. Zagier that for any even integer $k\geqslant 4$, we have the following restricted sum formulas
\begin{align*}
&\sum\limits_{i=2\atop i:\text{even}}^{k-1}\zeta_R(i,k-i)=\frac{3}{4}\zeta_R(k),\\
&\sum\limits_{i=2\atop i:\text{odd}}^{k-1}\zeta_R(i,k-i)=\frac{1}{4}\zeta_R(k),\\
&\sum\limits_{i=2\atop i:\text{even}}^{k-1}\zeta_R^{\star}(i,k-i)=\frac{2k-1}{4}\zeta_R(k),\\
&\sum\limits_{i=2\atop i:\text{odd}}^{k-1}\zeta_R^{\star}(i,k-i)=\frac{2k-3}{4}\zeta_R(k).
\end{align*}
We state a proof of the above equations. Using Euler's decomposition formulas
\begin{align*}
&z_r\shuffle z_s=\sum\limits_{i=r}^{r+s-1}\binom{i-1}{r-1}z_iz_{r+s-i}+\sum\limits_{i=s}^{r+s-1}\binom{i-1}{s-1}z_iz_{r+s-i},\\
&z_r\sshuffle z_s=\sum\limits_{i=r}^{r+s-1}\binom{i-1}{r-1}z_iz_{r+s-i}+\sum\limits_{i=s}^{r+s-1}\binom{i-1}{s-1}z_iz_{r+s-i}-\binom{r+s}{r}z_{r+s},
\end{align*}
one gets
\begin{align*}
&\sum\limits_{r=1}^{k-1}(-1)^r(z_r\shuffle z_{k-r}-z_r\ast z_{k-r})=-2\sum\limits_{i=2}^{k-1}(-1)^iz_iz_{k-i}+z_k,\\
&\sum\limits_{r=1}^{k-1}(-1)^r(z_r\sshuffle z_{k-r}-z_r\sast z_{k-r})=-2\sum\limits_{i=2}^{k-1}(-1)^iz_iz_{k-i}+z_k,
\end{align*}
which imply
\begin{align*}
&\sum\limits_{i=2}^{k-1}(-1)^i\zeta_R(i,k-i)=\frac{1}{2}\zeta_R(k),\\
&\sum\limits_{i=2}^{k-1}(-1)^i\zeta_R^{\star}(i,k-i)=\frac{1}{2}\zeta_R(k).
\end{align*}
Then using the sum formulas, one obtains the above restricted sum formulas.

There are also relations for $\zeta_R(r,s)$ and $\zeta_R^{\star}(r,s)$ with $r+s$ odd. For example, as shown in \cite{Gangl-Kaneko-Zagier}, for any odd integer $k\geqslant 3$ and integers $r,s$ with $r+s=k$ and $r,s\geqslant 2$, we have
\begin{align*}
\zeta_R(r,s)=&\frac{(-1)^r+1}{2}\zeta_R(r)\zeta_R(s)-(-1)^r\sum\limits_{r\leqslant i\leqslant k-2\atop i:\text{odd}}\binom{i-1}{r-1}\zeta_R(i)\zeta_R(k-i)\\
&-(-1)^r\sum\limits_{s\leqslant i\leqslant k-2\atop i:\text{odd}}\binom{i-1}{s-1}\zeta_R(i)\zeta_R(k-i)+\frac{1}{2}\left[(-1)^r\binom{k}{r}-1\right]\zeta_R(k)
\end{align*}
and
\begin{align*}
\zeta_R^{\star}(r,s)=&\frac{(-1)^r+1}{2}\zeta_R(r)\zeta_R(s)-(-1)^r\sum\limits_{r\leqslant i\leqslant k-2\atop i:\text{odd}}\binom{i-1}{r-1}\zeta_R(i)\zeta_R(k-i)\\
&-(-1)^r\sum\limits_{s\leqslant i\leqslant k-2\atop i:\text{odd}}\binom{i-1}{s-1}\zeta_R(i)\zeta_R(k-i)+\frac{1}{2}\left[(-1)^r\binom{k}{r}+1\right]\zeta_R(k).
\end{align*}
Moreover, since for odd integer $k\geqslant 3$, we have
\begin{align*}
2z_{k-1}z_1=&2(y\shuffle z_{k-1}-y\ast z_{k-1})+\sum\limits_{i=2}^{k-2}((-1)^iz_i\shuffle z_{k-i}\\
&-z_i\ast z_{k-i})+(k-1)z_k,
\end{align*}
we get
\begin{align*}
&\zeta_R(k-1,1)=\frac{k-1}{2}\zeta_R(k)-\sum\limits_{2\leqslant i\leqslant k-2\atop i:\text{odd}}\zeta_R(i)\zeta_R(k-i),\\
&\zeta_R^{\star}(k-1,1)=\frac{k+1}{2}\zeta_R(k)-\sum\limits_{2\leqslant i\leqslant k-2\atop i:\text{odd}}\zeta_R(i)\zeta_R(k-i).
\end{align*}


\subsection{Some evaluation formulas of multiple zeta value}

In this subsection, we show that some evaluation formulas, including the formulas for $\zeta_R(2k,\ldots,2k)$ and $\zeta^{\star}_R(2k,\ldots,2k)$, can be derived from the regularized double shuffle relations.

We take one $\lambda$ and then fix it, such that $\lambda^2=-24\zeta_R(2)$. Hence we have
$$\zeta_R(2)=-\frac{1}{24}\lambda^2.$$
Note that when $Z_R=Z$, we may take $\lambda=2\pi \sqrt{-1}$.

\begin{thm}\label{Thm:Ev-n-2}
Assume that the $\mathbb{K}$-linear map $Z_R:\mathfrak{h}^0\rightarrow R$ satisfies the regularized double shuffle relations, then for any nonnegative integer $n$, we have
\begin{align}
\zeta_{R}(\{2\}^n)=(-1)^n\frac{\lambda^{2n}}{4^n(2n+1)!}.
\label{Eq:Ev-n-2}
\end{align}
\end{thm}

To prove Theorem \ref{Thm:Ev-n-2}, we need the following result.

\begin{lem}\label{Lem:Insertion-4-2}
Assume that the $\mathbb{K}$-linear map $Z_R:\mathfrak{h}^0\rightarrow R$ satisfies the regularized double shuffle relations, then for any positive integer $n$, we have
\begin{align}
\sum\limits_{m=0}^{n-1}\zeta_{R}(\{2\}^m,4,\{2\}^{n-1-m})=\frac{2}{3}n(n+1)\zeta_{R}(\{2\}^{n+1}).
\label{Eq:Insertion-4-2}
\end{align}
\end{lem}

\noindent {\bf Proof of Theorem \ref{Thm:Ev-n-2}.} We proceed by induction on $n$. It is trivial to verify the result for $n=0$ and $n=1$.  Now assume that the formula \eqref{Eq:Ev-n-2} is valid for the positive integer $n$. Using the stuffle product, we have
$$\zeta_{R}(2)\zeta_{R}(\{2\}^n)=(n+1)\zeta_{R}(\{2\}^{n+1})+\sum\limits_{m=0}^{n-1}\zeta_{R}(\{2\}^m,4,\{2\}^{n-1-m}).$$
Hence by \eqref{Eq:Insertion-4-2}, we get
\begin{align*}
\zeta_{R}(2)\zeta_{R}(\{2\}^n)=\frac{(n+1)(2n+3)}{3}\zeta_{R}(\{2\}^{n+1}).
\end{align*}
Then using the induction assumption, we get the formula for $\zeta_{R}(\{2\}^{n+1})$, which finishes the proof. \qed

To prove Lemma \ref{Lem:Insertion-4-2}, we need the following two lemmas. The first one is immediately from the definitions of the stuffle products.

\begin{lem}\label{Lem:a-b-ast}
Let $a$ and $b$ be two positive integers. Then for any positive integer $n$, we have
$$\sum\limits_{m=0}^{n-1}(-1)^mz_{b+ma}\ast z_a^{n-1-m}=\sum\limits_{m=0}^{n-1}z_a^mz_bz_a^{n-1-m}$$
and
$$\sum\limits_{m=0}^{n-1}z_{b+ma}\sast z_a^{n-1-m}=\sum\limits_{m=0}^{n-1}z_a^mz_bz_a^{n-1-m}.$$
\end{lem}

The second one will need the restricted sum formula discussed in the last subsection.

\begin{lem}\label{Lem:2m-times-n-2}
Assume that the $\mathbb{K}$-linear map $Z_R:\mathfrak{h}^0\rightarrow R$ satisfies the regularized double shuffle relations. Then for any positive integer $n$, we have
\begin{align}
\sum\limits_{m=2}^{n+1}(-1)^m\zeta_{R}(2m)\zeta_{R}(\{2\}^{n+1-m})=\frac{2}{3}n(n+1)\zeta_{R}(\{2\}^{n+1}).
\label{Eq:Alternating-Sum-2m}
\end{align}
\end{lem}

\proof Let $t$ be a variable and set
$$F(t)=\exp\left(
\sum\limits_{n=1}^\infty \frac{(-1)^{n-1}}{n}\zeta_{R}(2n)t^n\right).$$
Then we get
\begin{align*}
&F'(t)=F(t)\cdot\left(\sum\limits_{n=1}^\infty (-1)^{n-1}\zeta_{R}(2n)t^{n-1}\right),\\
&F''(t)=F(t)\cdot\left\{\left(\sum\limits_{n=1}^\infty (-1)^{n-1}\zeta_{R}(2n)t^{n-1}\right)^2+\sum\limits_{n=1}^\infty (-1)^{n-1}(n-1)\zeta_{R}(2n)t^{n-2}\right\}.
\end{align*}
Using the stuffle products, we get
\begin{align*}
&\left(\sum\limits_{n=1}^\infty (-1)^{n-1}\zeta_{R}(2n)t^{n-1}\right)^2=\sum\limits_{n=2}^\infty(-1)^n\left(\sum\limits_{m+l=n\atop m,l\geqslant 1}\zeta_{R}(2m)\zeta_{R}(2l)\right)t^{n-2}\\
=&\sum\limits_{n=2}^\infty(-1)^n\sum\limits_{m+l=n\atop m,l\geqslant 1}\left(\zeta_{R}(2m,2l)+\zeta_{R}(2l,2m)+\zeta_{R}(2n)\right)t^{n-2}\\
=&2\sum\limits_{n=2}^\infty (-1)^n\left(\sum\limits_{m=1}^{n-1}\zeta_{R}(2m,2n-2m)\right)t^{n-2}+\sum\limits_{n=2}^\infty(-1)^n(n-1)\zeta_{R}(2n)t^{n-2}.
\end{align*}
Then by the restricted sum formula of double zeta values, we find
$$F''(t)=F(t)\cdot \frac{3}{2}\sum\limits_{n=2}^\infty (-1)^n\zeta_{R}(2n)t^{n-2}.$$

On the other hand, let $k=2$ in \eqref{Eq:kn-nk}, one has
$$F(t)=\sum\limits_{n=0}^\infty\zeta_R(\{2\}^n)t^n.$$
Therefore
$$F''(t)=\sum\limits_{n=0}^\infty n(n-1)\zeta_{R}(\{2\}^n)t^{n-2}.$$
Finally we get
\begin{align*}
&\left(\sum\limits_{n=2}^\infty (-1)^n\zeta_{R}(2n)t^{n-2}\right)\left(\sum\limits_{n=0}^\infty\zeta_{R}(\{2\}^n)t^n\right)=\frac{2}{3}\sum\limits_{n=0}^\infty n(n-1)\zeta_{R}(\{2\}^n)t^{n-2},
\end{align*}
which finishes the proof.
\qed

Now we come to prove Lemma \ref{Lem:Insertion-4-2}.

\noindent {\bf Proof of Lemma \ref{Lem:Insertion-4-2}.} Taking $a=2$ and $b=4$ in Lemma \ref{Lem:a-b-ast}, we obtain
$$\sum\limits_{m=0}^{n-1}(-1)^mz_{2(m+2)}\ast z_2^{n-1-m}=\sum\limits_{m=0}^{n-1}z_2^mz_4z_2^{n-1-m}.$$
Then \eqref{Eq:Insertion-4-2} follows from \eqref{Eq:Alternating-Sum-2m}.
\qed

Recall that
$$\Gamma_{R}(s)=\exp\left(\sum\limits_{n=2}^\infty\frac{(-1)^n}{n}\zeta_{R}(n)s^n\right)\in R[[s]]$$
is the gamma series associated to the map $Z_R:\mathfrak{h}^0\rightarrow R$. Then by Theorem \ref{Thm:Ev-n-2} and \cite[Theorem 4.5]{Li2010}, we get the following theorem.

\begin{thm}\label{Thm:Reflection-formula}
Assume that the $\mathbb{K}$-linear map $Z_R:\mathfrak{h}^0\rightarrow R$ satisfies the regularized double shuffle relations. Then the gamma series $\Gamma_{R}(s)$ satisfies the reflection formula
$$\Gamma_{R}(s)\Gamma_{R}(-s)=\frac{\lambda s}{e^{\lambda s/2}-e^{-\lambda s/2}}.$$
\end{thm}

Using \cite[Theorem 4.5]{Li2010}, we can obtain some other evaluation formulas. Before stating the formulas, we recall some notations.

Let $a,b,c$ be positive integers with $a+b=2c$. For any integers $m,n$ with $m\geqslant 2n\geqslant 0$, we denote by $I_{m,n}$ the set of all indices obtained from shuffling $(\{a,b\}^n)$ with $(\{c\}^{m-2n})$. For example, we have
$$I_{m,0}=\{(\{c\}^{m})\},\qquad I_{2n,n}=\{(\{a,b\}^{n})\}$$
and
$$I_{2n+1,n}=\{(c,\{a,b\}^{n}),(a,c,b,\{a,b\}^{n-1}),(a,b,c,\{a,b\}^{n-1}),\ldots,(\{a,b\}^{n},c)\}.$$
Define
\begin{align}
T_{m,n}=\sum\limits_{\mathbf{k}=(k_1,\ldots,k_m)\in I_{m,n}}z_{k_1}\ldots z_{k_m}\in\mathfrak{h}^1,
\label{Eq:Tmn}
\end{align}
then by \cite[Eq. (3.1)]{Yamamoto}, we have
\begin{align}
S(T_{m,n})=\sum\limits_{{2j+k+p=2n\atop i+k+l+p+q=m}\atop k,l,p,q\geqslant 0,i\geqslant 2j\geqslant 0}(-1)^{i+k}
\binom{k+l}{k}\binom{p+q}{p}T_{i,j}\ast S(z_c^{k+l})\ast S(z_c^{p+q}).
\label{Eq:Insertion-sum-h1}
\end{align}
Finally, let $\{B_n\}$ be the Bernoulli numbers defined by
$$\sum\limits_{n=0}^\infty \frac{B_n}{n!}t^n=\frac{t}{e^t-1}.$$

\begin{thm}\label{Thm:Ev-2}
Assume that $Z_R:\mathfrak{h}^0\rightarrow R$ satisfies the regularized double shuffle relations, then
\begin{itemize}
  \item[(1)] for any nonnegative integer $n$, we have
  $$\zeta_{R}(2n)=-\frac{\lambda^{2n}}{2(2n)!}B_{2n},$$
  where we set $\zeta_R(0)=-\frac{1}{2}$;
    \item[(2)] for any positive integer $n$, we have $$\zeta_{R}^{\star}(\{2\}^n)=(2^{1-2n}-1)\frac{\lambda^{2n}}{(2n)!}B_{2n}=2(1-2^{1-2n})\zeta_{R}(2n);$$
  \item[(3)] for any positive integer $k$ and any nonnegative integer $n$, we have
      $$\zeta_{R}(\{2k\}^n)=C_{n}^{(k)}\frac{\lambda^{2nk}}{(2nk)!},$$
      where $C_n^{(k)}\in\mathbb{Q}$ are given by $C_0^{(k)}=1$ and the recursive formula
$$C_{n}^{(k)}=\frac{1}{2n}\sum\limits_{m=1}^n(-1)^m\binom{2nk}{2mk}B_{2mk}C_{n-m}^{(k)},\quad (\forall n\geqslant 1);$$
\item[(4)] for any positive integer $k$ and any nonnegative integer $n$, we have
      $$\zeta_{R}^{\star}(\{2k\}^n)=C_{n}^{\star,(k)}\frac{\lambda^{2nk}}{(2nk)!},$$
      where $C_n^{\star,(k)}\in\mathbb{Q}$ are given by $C_0^{\star,(k)}=1$ and the recursive formula
$$C_{n}^{\star,(k)}=-\frac{1}{2n}\sum\limits_{m=1}^n\binom{2nk}{2mk}B_{2mk}C_{n-m}^{\star,(k)},\quad (\forall n\geqslant 1);$$
  \item[(5)] for any integers $m,n$ with $m\geqslant 2n\geqslant 0$, we have
  \begin{align*}
 &\sum\limits_{m_0+m_1+\cdots+m_{2n}=m-2n\atop m_0,\ldots,m_{2n}\geqslant 0}\zeta_{R}(\{2\}^{m_0},3,\{2\}^{m_1},1,\{2\}^{m_2},3,\{2\}^{m_3},1
     ,\ldots,1,\\
     &\qquad\qquad\{2\}^{m_{2n-2}},3,\{2\}^{m_{2n-1}},1,\{2\}^{m_{2n}})=(-1)^m\frac{2}{4^m(2m+2)!}
     \binom{m+1}{2n+1}\lambda^{2m}.
     \end{align*}
In particular, for any nonnegative integer $n$, we have
$$\zeta_{R}(\{3,1\}^{n})=\frac{\lambda^{4n}}{2^{4n-1}(4n+2)!}$$
and
\begin{align*}
&\zeta_{R}(2,\{3,1\}^{n})+\zeta_{R}(3,2,1,\{3,1\}^{n-1})
+\zeta_{R}(3,1,2,\{3,1\}^{n-1})\\
&+\cdots+\zeta_{R}(\{3,1\}^{n},2)=-\frac{\lambda^{4n+2}}{4^{2n+1}(4n+3)!};
\end{align*}
  \item[(6)] for any integers $m,n$ with $m\geqslant 2n\geqslant 0$, we have
\begin{align*}
 &\sum\limits_{m_0+m_1+\cdots+m_{2n}=m-2n\atop m_0,\ldots,m_{2n}\geqslant 0}\zeta_{R}^{\star}(\{2\}^{m_0},3,\{2\}^{m_1},1,\{2\}^{m_2},3,\{2\}^{m_3},1
     ,\ldots,1,\\
     &\qquad\qquad\{2\}^{m_{2n-2}},3,\{2\}^{m_{2n-1}},1,\{2\}^{m_{2n}})\\
=&\sum\limits_{{2j+k+p=2n\atop i+k+l+p+q=m}\atop k,l,p,q\geqslant 0,i\geqslant 2j\geqslant 0}(-1)^{k}
\binom{k+l}{k}\binom{p+q}{p}\binom{i+1}{2j+1}\frac{\beta_{k+l}\beta_{p+q}}{2^{2i-1}(2i+2)!}\lambda^{2m},
\end{align*}
where $\beta_k=\frac{2^{1-2k}-1}{(2k)!}B_{2k}$. In particular, for any nonnegative integer $n$, we have
$$\zeta_{R}^{\star}(\{3,1\}^{n})=\sum\limits_{2j+k+p=2n\atop j,k,p\geqslant 0}
(-1)^k\frac{\beta_k\beta_p}{2^{4j-1}(4j+2)!}\lambda^{4n}$$
and
\begin{align*}
&\zeta_{R}^{\star}(2,\{3,1\}^{n})+\zeta_{R}^{\star}(3,2,1,\{3,1\}^{n-1})
+\zeta_{R}^{\star}(3,1,2,\{3,1\}^{n-1})\\
&+\cdots+\zeta_{R}^{\star}(\{3,1\}^{n},2)=\sum\limits_{j=0}^n
\left\{\frac{1}{4^{2j+1}(4j+3)!}\sum\limits_{k+p=2(n-j)\atop k,p\geqslant 0}(-1)^k\beta_k\beta_p\right.\\
&\qquad\qquad\left.-\frac{1}{4^{2j-1}(4j+2)!}\sum\limits_{k+p=2(n-j)+1\atop k\geqslant 1,p\geqslant 0}(-1)^kk\beta_k\beta_p\right\}
\lambda^{4n+2}.
\end{align*}
\end{itemize}
\end{thm}

\proof We get (1), (3) and (5) from Theorem \ref{Thm:Ev-n-2} and \cite[Theorem 4.5]{Li2010}. By \eqref{Eq:MZVk-MZSVk} and Theorem \ref{Thm:Reflection-formula}, we get
\begin{align}
\sum\limits_{n=0}^\infty\zeta_R^{\star}(\{2\}^n)s^{2n}=&\Gamma_{R}(s)\Gamma_{R}(-s)=\frac{\lambda s}{e^{\lambda s/2}-e^{-\lambda s/2}}
\label{Eq:Generating-2-n}\\
=&2\frac{\lambda s/2}{e^{\lambda s/2}-1}-\frac{\lambda s}{e^{\lambda s}-1},\nonumber
\end{align}
which deduces (2).

Taking $a=b=2k$ in Lemma \ref{Lem:a-b-ast}, we get
$$\sum\limits_{m=1}^nz_{2mk}\sast z_{2k}^{n-m}=nz_{2k}^n,$$
which implies
$$\zeta_R^{\star}(\{2k\}^n)=\frac{1}{n}\sum\limits_{m=1}^n\zeta_R(2mk)\zeta_R^{\star}(\{2k\}^{n-m}).$$
Hence by (1), we have
$$\zeta_R^{\star}(\{2k\}^n)
=-\frac{1}{2n}\sum\limits_{m=1}^n\frac{\lambda^{2mk}}{(2mk)!}B_{2mk}\zeta_R^{\star}(\{2k\}^{n-m}).$$
Then we get (4) by induction on $n$.

Taking $a=3,b=1,c=2$ in \eqref{Eq:Insertion-sum-h1} and using (2) and (5), we get (6).
\qed

Let $n=0$ in the item (6) of Theorem \ref{Thm:Ev-2}, we get
$$\zeta_{R}^{\star}(\{2\}^m)=\sum\limits_{i+l+q=m\atop i,l,q\geqslant 0}\frac{(2^{1-2l}-1)(2^{1-2q}-1)}{4^i(2i+1)!(2l)!(2q)!}B_{2l}B_{2q}\lambda^{2m}.$$
Comparing with the item (2) of Theorem \ref{Thm:Ev-2}, we obtain a relation among Bernoulli numbers.

\begin{cor}
For any nonnegative integer $n$, we have
$$
\sum\limits_{i+j+k=n\atop i,j,k\geqslant 0}\frac{(2^{1-2j}-1)(2^{1-2k}-1)}{4^i(2i+1)!(2j)!(2k)!}B_{2j}B_{2k}=\frac{2^{1-2n}-1}{(2n)!}B_{2n}.
$$
\end{cor}

In fact, one can use
\begin{align*}
&\sum\limits_{n=0}^\infty\frac{2^{1-2n}-1}{(2n)!}B_{2n}t^{2n}=\frac{t}{e^{t/2}-e^{-t/2}},\\
&\sum\limits_{n=0}^\infty\frac{1}{4^n(2n+1)!}t^{2n}=\frac{e^{t/2}-e^{-t/2}}{t}
\end{align*}
to prove the above corollary.

If $Z_R=Z$ and $\lambda=2\pi\sqrt{-1}$, we get the evaluation formulas for multiple zeta values and multiple zeta-star values:
\begin{itemize}
  \item[(i)] the item (1) of Theorem \ref{Thm:Ev-2} becomes Euler's formula
  $$\zeta(2n)=(-1)^{n+1}\frac{2^{2n-1}\pi^{2n}}{(2n)!}B_{2n};$$
  \item[(ii)] Theorem \ref{Thm:Ev-n-2} is (\cite{Borwein-Bradley-Broadhurst,Hoffman1997})
  $$\zeta(\{2\}^{n})=\frac{\pi^{2n}}{(2n+1)!};$$
  \item[(iii)] the item (2) of Theorem \ref{Thm:Ev-2} is (\cite{Muneta2008,Zlobin})
  $$\zeta^{\star}(\{2\}^n)=(-1)^n\frac{(2-2^{2n})\pi^{2n}}{(2n)!}B_{2n}=2(1-2^{1-2n})\zeta(2n);$$
  \item[(iv)] the item (3) of Theorem \ref{Thm:Ev-2} becomes a formula of $\zeta(\{2k\}^n)$ (\cite[Section 3]{Borwein-Bradley-Broadhurst});
  \item[(v)] the item (4) of Theorem \ref{Thm:Ev-2} becomes a formula of $\zeta^{\star}(\{2k\}^n)$;
  \item[(vi)] the item (5) of Theorem \ref{Thm:Ev-2} contains the formula
  $$\zeta(T_{m,n})=\frac{2\pi^{2m}}{(2m+2)!}\binom{m+1}{2n+1},$$
  and formulas of $\zeta(\{3,1\}^{n})$ and $\zeta(T_{2n+1,n})$ with $a=3,b=1,c=2$ (\cite{Borwein-Bradley-Broadhurst}, \cite[Theorem 1 and Theorem 2]{Borwein-Bradley-Broadhurst-Lisonek}, \cite[Corollary 5.1]{Bowman-Bradley});
  \item[(vii)] the item (6) of Theorem \ref{Thm:Ev-2} contains the formulas of $\zeta^{\star}(T_{m,n})$, $\zeta^{\star}(\{3,1\}^{n})$ and $\zeta^{\star}(T_{2n+1,n})$ with $a=3,b=1,c=2$ (\cite[Theorem B]{Muneta2008},\cite[Theorem 1.1]{Yamamoto}).
\end{itemize}
Hence the above evaluation formulas for multiple zeta values and multiple zeta-star values can be derived from the regularized double shuffle relations.

The items (3) and (4) of Theorem \ref{Thm:Ev-2} give the formulas of  $\zeta_R(\{2k\}^n)$ and $\zeta_R^{\star}(\{2k\}^n)$. While there are constants defined by recursive relations in the formulas. In fact, we have

\begin{thm}\label{Thm:Ev-2k}
Assume that the $\mathbb{K}$-linear map $Z_R:\mathfrak{h}^0\rightarrow R$ satisfies the regularized double shuffle relations. Then for any positive integer $k$ and any nonnegative integer $n$, we have
\begin{align}
&\zeta_{R}(\{2k\}^n)=(-1)^n\sum\limits_{n_0+\cdots+n_{k-1}=nk\atop n_0,\ldots,n_{k-1}\geqslant 0}
\frac{\rho_k^{\sum\limits_{j=0}^{k-1}2jn_j}}{(2n_0+1)!\cdots (2n_{k-1}+1)!}\frac{\lambda^{2nk}}{4^{nk}},
\label{Eq:z-2k-n}\\
&\zeta_{R}^{\star}(\{2k\}^n)=\sum\limits_{n_0+\cdots+n_{k-1}=nk\atop n_0,\ldots,n_{k-1}\geqslant 0}
\left(\prod\limits_{j=0}^{k-1}\frac{B_{2n_j}}{(2n_j)!}(2-4^{n_j})\right)\rho_k^{\sum\limits_{j=0}^{k-1}2jn_j}\frac{\lambda^{2nk}}{4^{nk}},
\label{Eq:zs-2k-n}
\end{align}
where $\rho_{k}=e^{\frac{\pi\sqrt{-1}}{k}}$.
\end{thm}

Note that as \cite[Theorem A]{Muneta2008}, one can show that the coefficients of $\frac{\lambda^{2nk}}{4^{nk}}$ in the right-hand sides of \eqref{Eq:z-2k-n} and \eqref{Eq:zs-2k-n} are in fact rational numbers. And if $Z_R=Z$ and $\lambda=2\pi \sqrt{-1}$, we have (\cite[Theorem A]{Muneta2008})
$$\zeta^{\star}(\{2k\}^n)=(-1)^{nk}\sum\limits_{n_0+\cdots+n_{k-1}=nk\atop n_0,\ldots,n_{k-1}\geqslant 0}
\left(\prod\limits_{j=0}^{k-1}\frac{B_{2n_j}}{(2n_j)!}(2-4^{n_j})\right)\rho_k^{\sum\limits_{j=0}^{k-1}2jn_j}\pi^{2nk}.$$

Let $k=2$. Then $\rho_k=\sqrt{-1}$. Hence from \eqref{Eq:z-2k-n} and \eqref{Eq:zs-2k-n}, we get
\begin{align}
&\zeta_{R}(\{4\}^n)=\sum\limits_{i=0}^{2n}\frac{(-1)^{n+i}}{(2i+1)!(4n-2i+1)!}\frac{\lambda^{4n}}{4^{2n}},
\label{Eq:z-4-n}\\
&\zeta_{R}^{\star}(\{4\}^n)=\sum\limits_{i=0}^{2n}(-1)^{i}\frac{(2-4^{i})(2-4^{2n-i})B_{2i}B_{4n-2i}}{(2i)!(4n-2i)!}
\frac{\lambda^{4n}}{4^{2n}}.
\label{Eq:zs-4-n}
\end{align}
While since
$$\sum\limits_{i=0}^{2n}(-1)^{i}\binom{4n+2}{2i+1}=(-1)^n2^{2n+1},$$
\eqref{Eq:z-4-n} reduces to
\begin{align}
\zeta_{R}(\{4\}^n)=\frac{2\lambda^{4n}}{4^{n}(4n+2)!}.
\label{Eq:z-4-n-new}
\end{align}
If $Z_R=Z$ and $\lambda=2\pi\sqrt{-1}$, we get (\cite[p. 7]{Borwein-Bradley-Broadhurst})
$$\zeta(\{4\}^n)=\frac{2^{2n+1}\pi^{4n}}{(4n+2)!}.$$

Now we come to give a proof of Theorem \ref{Thm:Ev-2k}. Let $n$ be any nonnegative integer. Then we have
\begin{align}
\sum\limits_{j=0}^{k-1}\rho_{k}^{2nj}=
\begin{cases}
0 & \text{if\;}k\nmid n,\\
k & \text{if\;}k\mid n.
\end{cases}
\label{Eq:sum-rho-k}
\end{align}

\begin{lem}\label{Lem:sum-Bernoulli}
Let $t$ be a variable. Then we have
$$\sum\limits_{m=0}^{\infty}\frac{\lambda^{2mk}}{(2mk)!}B_{2mk}t^{2mk}
=\frac{1}{k}\sum\limits_{j=0}^{k-1}\left(\frac{\lambda\rho_{k}^{j}t}{2}
\frac{e^{\lambda\rho_{k}^{j}t}+1}{e^{\lambda\rho_{k}^{j}t}-1}\right).$$
\end{lem}

\proof Since
$$\sum\limits_{n=0}^\infty\frac{B_{2n}}{(2n)!}t^{2n}=\frac{t}{e^t-1}+\frac{t}{2}=
\frac{t}{2}\frac{e^t+1}{e^t-1},$$
we have
$$\sum\limits_{j=0}^{k-1}\left(\frac{\lambda\rho_{k}^{j}t}{2}
\frac{e^{\lambda\rho_{k}^{j}t}+1}{e^{\lambda\rho_{k}^{j}t}-1}\right)
=\sum\limits_{n=0}^{\infty}\frac{\lambda^{2n}}{(2n)!}B_{2n}t^{2n}
\sum\limits_{j=0}^{k-1}\rho_{k}^{2nj}.$$
Then we get the result from \eqref{Eq:sum-rho-k}.
\qed

Now we can compute the generating functions of $\zeta_R(\{2k\}^n)$ and $\zeta_R^{\star}(\{2k\}^n)$.

\begin{prop}\label{Prop:generating-z-2k-n}
Assume that the $\mathbb{K}$-linear map $Z_R:\mathfrak{h}^0\rightarrow R$ satisfies the regularized double shuffle relations. Let $k$ be a positive integer and $t$ be a variable. Then we have
\begin{align}
\sum\limits_{n=0}^{\infty}(-1)^n\zeta_{R}(\{2k\}^n)t^{2nk}=
\prod\limits_{j=0}^{k-1}\frac{e^{\frac{1}{2}\lambda\rho_{k}^jt}
-e^{-\frac{1}{2}\lambda\rho_{k}^jt}}{\lambda\rho_{k}^jt}
\label{Eq:Generating-z-2k}
\end{align}
and
\begin{align}
\sum\limits_{n=0}^{\infty}\zeta_{R}^{\star}(\{2k\}^n)t^{2nk}=
\prod\limits_{j=0}^{k-1}\frac{\lambda\rho_{k}^jt}{e^{\frac{1}{2}\lambda\rho_{k}^jt}
-e^{-\frac{1}{2}\lambda\rho_{k}^jt}}.
\label{Eq:Generating-zs-2k}
\end{align}
\end{prop}

\proof
Taking $a=b=2k$ in Lemma \ref{Lem:a-b-ast}, we have
\begin{align*}
\sum\limits_{m=1}^{n}(-1)^{m-1}z_{2mk}\ast z_{2k}^{n-m}=nz_{2k}^{n},
\end{align*}
which induces
$$\sum\limits_{m=1}^{n}(-1)^{m-1}\zeta_{R}(2mk)\zeta_{R}(\{2k\}^{n-m})=n\zeta_{R}(\{2k\}^n).$$
By item (1) of Theorem \ref{Thm:Ev-2}, we obtain
$$2n\zeta_{R}(\{2k\}^n)
=\sum\limits_{m=1}^n(-1)^{m}\frac{\lambda^{2mk}}{(2mk)!}B_{2mk}\zeta_{R}(\{2k\}^{n-m}),$$
where $n$ is any nonnegative integer. Therefore we have
\begin{align*}
&\sum\limits_{n=0}^{\infty}(-1)^n(2n+1)\zeta_{R}(\{2k\}^n)t^{2nk}\\
&=\sum\limits_{n=0}^{\infty}(-1)^nt^{2nk}\sum\limits_{m=0}^n(-1)^{m}
\frac{\lambda^{2mk}}{(2mk)!}B_{2mk}\zeta_{R}(\{2k\}^{n-m})\\
&=\sum\limits_{m=0}^{\infty}\frac{\lambda^{2mk}}{(2mk)!}B_{2mk}t^{2mk}
\sum\limits_{n=0}^{\infty}(-1)^n\zeta_{R}(\{2k\}^n)t^{2nk}.
\end{align*}

Set
$$f(t)=\sum\limits_{n=0}^{\infty}(-1)^n\zeta_{R}(\{2k\}^n)t^{2nk}.$$
By Lemma \ref{Lem:sum-Bernoulli}, we get
\begin{align*}
\sum\limits_{j=0}^{k-1}\left(\frac{\lambda\rho_{k}^{j}t}{2}
\frac{e^{\lambda\rho_{k}^{j}t}+1}{e^{\lambda\rho_{k}^{j}t}-1}\right)f(t)&=
\sum\limits_{n=0}^{\infty}(-1)^n(2n+1)k\zeta_{R}(\{2k\}^n)t^{2nk}\\
&=\frac{(t^kf(t))^{'}}{t^{k-1}}=kf(t)+tf'(t).
\end{align*}
Hence $f(t)$ satisfies the differential equation
$$f'(t)=\sum\limits_{j=0}^{k-1}\left(\frac{\lambda\rho_{k}^{j}}{2}
\frac{e^{\lambda\rho_{k}^{j}t}+1}{e^{\lambda\rho_{k}^{j}t}-1}-\frac{1}{t}\right)f(t).$$
Solving the differential equation, we get
$$f(t)=C\prod\limits_{j=0}^{k-1}\frac{e^{\frac{1}{2}\lambda\rho_{k}^jt}-
e^{-\frac{1}{2}\lambda\rho_{k}^jt}}{t},$$
where $C$ is a constant. Since $f(0)=1$, we get
$$C=\left(\prod\limits_{j=0}^{k-1}(\lambda\rho_k^j)\right)^{-1},$$
which induces \eqref{Eq:Generating-z-2k}.

We get \eqref{Eq:Generating-zs-2k} from \eqref{Eq:MZVk-MZSVk} and \eqref{Eq:Generating-z-2k}.
\qed

Note that if $Z_R=Z$ and $\lambda=2\pi\sqrt{-1}$, \eqref{Eq:Generating-z-2k} is equivalent to \cite[Eq. (34)]{Borwein-Bradley-Broadhurst}.

Now we give a proof of Theorem \ref{Thm:Ev-2k}.

\noindent {\bf Proof of Theorem \ref{Thm:Ev-2k}.} We get the result from the expansions
\begin{align*}
&\frac{e^{\frac{s}{2}}-e^{-\frac{s}{2}}}{s}=\sum\limits_{n=0}^{\infty}\frac{s^{2n}}{4^n(2n+1)!},\\
&\frac{s}{e^{\frac{s}{2}}-e^{-\frac{s}{2}}}=\sum\limits_{n=0}^\infty\frac{B_{2n}}{4^n(2n)!}\left(2-4^{n}\right)s^{2n},
\end{align*}
and Proposition \ref{Prop:generating-z-2k-n}.
\qed


\subsection{Ohno-Zagier relation and some corollaries}

In \cite{Li2013}, the first named author proved that the Ohno-Zagier relation (\cite{Ohno-Zagier}) can be deduced from the regularized double shuffle relations. In fact, set
$$X_0(k,n,h)=\sum\limits_{\wt(\mathbf{k})=k,\dep(\mathbf{k})=n,\height(\mathbf{k})=h\atop \mathbf{k}:\text{admissible}}\zeta_{R}(\mathbf{k}).$$
Then we have

\begin{thm}[\cite{Li2013}]\label{Thm:EDS-OZ}
Assume that $Z_R:\mathfrak{h}^0\rightarrow R$ satisfies the regularized double shuffle relations. Let $u,s,t$ be variables. Then we have
\begin{align}
&\sum\limits_{k\geqslant n+h\atop n\geqslant h\geqslant 1}X_0(k,n,h)u^{k-n-h}s^{n-h}t^{h-1}\nonumber\\
=&\frac{1}{us-t}\left\{1-\exp\left(\sum\limits_{n=2}^\infty \frac{\zeta_{R}(n)}{n} (u^n+s^n-\alpha^n-\beta^n)\right)\right\}.
\label{Eq:Ohno-Zagier}
\end{align}
Here $\alpha$ and $\beta$ are determined by $\alpha+\beta=u+s$ and $\alpha\beta=t$.
\end{thm}

Denote the left-hand side of \eqref{Eq:Ohno-Zagier} by $G_0(u,s,t)$. Since the right-hand side of \eqref{Eq:Ohno-Zagier} is symmetric about $u$ and $s$, we get

\begin{cor}[\cite{Kajikawa}]\label{Cor:duality}
Assume that the $\mathbb{K}$-linear map $Z_R:\mathfrak{h}^0\rightarrow R$ satisfies the regularized double shuffle relations. Then for any integers $k,n,h$ with $k\geqslant n+h$ and $n\geqslant h\geqslant 1$, we have
$$\text{\lq\lq duality of \;} X_0(k,n,h)\text{\rq\rq}=X_0(k,n,h).$$
\end{cor}

In above, for any $w\in\mathfrak{h}^0$, the duality of $Z_R(w)$ is $Z_R(\tau(w))$.

As in \cite{Ohno-Zagier}, we have some consequences of Theorem \ref{Thm:EDS-OZ}.
Let $u=s=0$ in Theorem \ref{Thm:EDS-OZ}. Since $\alpha,\beta=\pm\sqrt{-t}$ in this case, we get
$$\sum\limits_{n=0}^\infty\zeta_{R}(\{2\}^n)t^n=\exp\left(\sum\limits_{n=1}^\infty \frac{(-1)^{n-1}}{n}\zeta_{R}(2n)t^n\right),$$
which is of course a special case of \eqref{Eq:kn-nk}.

Let $t=0$ in Theorem \ref{Thm:EDS-OZ}. We take $\alpha=u+s$ and $\beta=0$. Then we have the Aomoto-Drinfel'd-Zagier relation
$$\sum\limits_{m,n=1}^\infty\zeta_{R}(m+1,\{1\}^{n-1})u^ms^n=1-
\exp\left(\sum\limits_{n=2}^\infty\frac{\zeta_{R}(n)}{n}(u^n+s^n-(u+s)^n)\right).$$

Finally, let $u=-s$ in Theorem \ref{Thm:EDS-OZ}. Then we have
$$G_0(u,-u,t)=\sum\limits_{k\geqslant 2h\geqslant 2}(-1)^h\left(\sum\limits_{\wt(\mathbf{k})=k,\height(\mathbf{k})=h\atop \mathbf{k}:\text{admissible}}(-1)^{\dep(\mathbf{k})}\zeta_{R}(\mathbf{k})\right)u^{k-2h}t^{h-1}.$$
On the other hand, since $\alpha,\beta=\pm\sqrt{-t}$, we have
\begin{align*}
&G_0(u,-u,t)=\frac{1}{u^2+t}\left\{\exp\left(\sum\limits_{n=1}^\infty \frac{\zeta_{R}(2n)}{n}(u^{2n}-(-t)^n)\right)-1\right\}\\
=&\exp\left(\sum\limits_{n=1}^\infty \frac{\zeta_{R}(2n)}{n}u^{2n}\right)\frac{\exp\left(-\sum\limits_{n=1}^\infty \frac{\zeta_{R}(2n)}{n}(-t)^{n}\right)
-\exp\left(-\sum\limits_{n=1}^\infty \frac{\zeta_{R}(2n)}{n}u^{2n}\right)}{u^2+t}\\
=&\Gamma_{R}(u)\Gamma_{R}(-u)\frac{1}{u^2+t}\left(\frac{1}{\Gamma_{R}((-t)^{1/2})\Gamma_{R}(-(-t)^{1/2})}
-\frac{1}{\Gamma_{R}(u)\Gamma_{R}(-u)}\right).
\end{align*}

By Theorem \ref{Thm:Reflection-formula}, we have the expansions
\begin{align}
\Gamma_{R}(s)\Gamma_{R}(-s)=&\sum\limits_{n=0}^\infty \frac{B_{2n}}{(2n)!}(2^{1-2n}-1)\lambda^{2n}s^{2n},
\label{Eq:gamma-extension}
\end{align}
and
\begin{align}
\frac{1}{\Gamma_{R}(s)\Gamma_{R}(-s)}=&\sum\limits_{n=0}^\infty \frac{\lambda^{2n}}{2^{2n}(2n+1)!}s^{2n}.
\label{Eq:gamma-inver-extension}
\end{align}
By \eqref{Eq:gamma-inver-extension}, we have
\begin{align*}
&\frac{1}{u^2+t}\left(\frac{1}{\Gamma_{R}((-t)^{1/2})\Gamma_{R}(-(-t)^{1/2})}-\frac{1}{\Gamma_{R}(u)\Gamma_{R}(-u)}\right)\\
=&\frac{1}{u^2+t}\sum\limits_{n=0}^\infty \frac{\lambda^{2n}}{2^{2n}(2n+1)!}((-t)^n-u^{2n})\\
=&-\sum\limits_{n=1}^\infty \frac{\lambda^{2n}}{2^{2n}(2n+1)!}\sum\limits_{m=1}^n u^{2n-2m}(-t)^{m-1}\\
=&\sum\limits_{n\geqslant m\geqslant 1}\frac{\lambda^{2n}}{2^{2n}(2n+1)!}(-1)^mu^{2n-2m}t^{m-1}.
\end{align*}
Hence using \eqref{Eq:gamma-extension}, we find $G_0(u,-u,t)$ is
$$\sum\limits_{k\geqslant h\geqslant 1}(-1)^h\frac{\lambda^{2k}}{2^{2k}(2k+1)!}\left(\sum\limits_{n=0}^{k-h}\binom{2k+1}{2n}(2-2^{2n})B_{2n}\right)u^{2k-2h}t^{h-1}.$$
Therefore, for any integers $k,h$ with $k\geqslant h\geqslant 1$, we have
$$\sum\limits_{\wt(\mathbf{k})=2k+1,\height(\mathbf{k})=h\atop \mathbf{k}:\text{admissible}}(-1)^{\dep(\mathbf{k})}\zeta_{R}(\mathbf{k})=0$$
and the Le-Murakami relation
$$\sum\limits_{\wt(\mathbf{k})=2k,\height(\mathbf{k})=h\atop \mathbf{k}:\text{admissible}}(-1)^{\dep(\mathbf{k})}\zeta_{R}(\mathbf{k})=\frac{\lambda^{2k}}{2^{2k}(2k+1)!}\sum\limits_{n=0}^{k-h}
\binom{2k+1}{2n}(2-2^{2n})B_{2n}.$$

Note that the Le-Murakami relation for multiple zeta values was first proved in \cite{Le-Murakami}.


\subsection{Restricted sum formulas of M. E. Hoffman and generalizations}

In \cite{Hoffman2017}, M. E. Hoffman considered the restricted sum
$$\sum\limits_{k_1+\cdots+k_n=k\atop k_i\geqslant 1}\zeta(2k_1,\ldots,2k_n),$$
and used the symmetric functions to give two formulas of this restricted sum. We find that one can translate the formulas in \cite{Hoffman2017} deduced from the properties of symmetric functions  to identities in the algebra $\mathfrak{h}$. Hence it is natural to consider more general sums
$$\sum\limits_{k_1+\cdots+k_n=k\atop k_i\geqslant 1}\zeta(ak_1,\ldots,ak_n), \quad \sum\limits_{k_1+\cdots+k_n=k\atop k_i\geqslant 1}\zeta^{\star}(ak_1,\ldots,ak_n),$$
where $a>1$ is an integer.

Now let $a$ be a fixed positive integer. We set
\begin{align*}
&E(t)=\sum\limits_{j=0}^\infty z_a^jt^j\in\mathfrak{h}^1[[t]],\\
&H(t)=\sum\limits_{j=0}^{\infty} S(z_a^j)t^j=S(E(t))\in\mathfrak{h}^1[[t]],\\
&\overline{H}(t)=\sum\limits_{j=0}^\infty S^{-1}(z_a^j)t^j=S^{-1}(E(t))\in\mathfrak{h}^1[[t]],\\
&P(t)=\sum\limits_{j=1}^\infty z_{aj}t^{j-1}\in\mathfrak{h}^1[[t]].
\end{align*}

\begin{lem}\label{Lem:RelSymFunc}
We have
\begin{itemize}
  \item[(1)] $E(-t)\ast H(t)=1$, $\overline{H}(-t)\sast E(t)=1$;
  \item[(2)] $P(t)=-H(t)\ast \frac{d}{dt}E(-t)=E(-t)\ast \frac{d}{dt}H(t)$;
  \item[(3)] $P(t)=-E(t)\sast \frac{d}{dt}\overline{H}(-t)=\overline{H}(-t)\sast \frac{d}{dt}E(t)$;
  \item[(4)] $\frac{d}{dt}E(-t)=-E(-t)\ast P(t)$, $\frac{d}{dt}H(t)=H(t)\ast P(t)$;
  \item[(5)] $\frac{d}{dt}E(t)=E(t)\sast P(t)$, $\frac{d}{dt}\overline{H}(-t)=-\overline{H}(-t)\sast P(t)$.
\end{itemize}
\end{lem}

\proof From \eqref{Eq:k-ast-Sk}, we get (1). By \eqref{Eq:Exp-ast}, we have
$$\exp_{\ast}\left(-\sum\limits_{n=1}^\infty \frac{z_{na}}{n}t^n\right)=E(-t).$$
Hence we get
$$E(-t)\ast \left(-\sum\limits_{n=1}^\infty z_{na}t^{n-1}\right)=\frac{d}{dt}E(-t),$$
which is
$$-E(-t)\ast P(t)=\frac{d}{dt}E(-t).$$
By (1), we have
$$\frac{d}{dt}E(-t)\ast H(t)=-E(-t)\ast\frac{d}{dt}H(t).$$
Then we get (2). Applying the map $S^{-1}$, we get (3) from (2). Finally, we get (4) from (1), (2), and get (5) from (1), (3).
\qed

Now for integers $k,n$ with $k\geqslant n\geqslant 1$, we define
$$N_{k,n}=\sum\limits_{k_1+\cdots+k_n=k\atop k_i\geqslant 1}z_{ak_1}\cdots z_{ak_n}\in\mathfrak{h}^1
$$
and their generating function
$$\mathcal{F}(t,s)=1+\sum\limits_{k\geqslant n\geqslant 1}N_{k,n}t^ks^n\in\mathfrak{h}^1[[t,s]],$$
where $t,s$ are variables.

Corresponding to \cite[Lemma 2.1]{Hoffman2017}, we have the following key result.

\begin{prop}\label{Prop:ResSum-E-H}
We have
\begin{align}
&E((s-1)t)=E(-t)\ast\mathcal{F}(t,s),
\label{Eq:E-ast-F}\\
&E((s+1)t)=E(t)\sast\mathcal{F}(t,s).
\label{Eq:E-sast-F}
\end{align}
Moreover, we have
\begin{align}
\mathcal{F}(t,s)=E((s-1)t)\ast H(t)=E((s+1)t)\sast \overline{H}(-t).
\label{Eq:F-ast-E}
\end{align}
\end{prop}

In fact, the two expressions in \eqref{Eq:F-ast-E} for $\mathcal{F}(t,s)$ are equivalent.  To see this, we need the following lemma (compare with \cite[Lemma 2.3]{Hoffman2017} and \cite[Lemma 4.1]{Gencev}).

\begin{lem}\label{Lem:S-F}
For any integers $k,n$ with $k\geqslant n\geqslant 1$, we have
\begin{align}
S^{-1}(N_{k,n})=\sum\limits_{i=1}^n(-1)^{n-i}\binom{k-i}{k-n}N_{k,i}
\label{Eq:S-Inverse-N}
\end{align}
and
\begin{align}
S(N_{k,n})=\sum\limits_{i=1}^n\binom{k-i}{k-n}N_{k,i},
\label{Eq:S-N}
\end{align}
which in terms of generating functions are 
\begin{align}
S^{-1}(\mathcal{F}(t,s))=\mathcal{F}\left(t(1-s),\frac{s}{1-s}\right)
\label{Eq:Generating-S-Inverse-N}
\end{align}
and
\begin{align}
S(\mathcal{F}(t,s))=\mathcal{F}\left(t(s+1),\frac{s}{s+1}\right).
\label{Eq:Generating-S-N}
\end{align}
\end{lem}

\proof We prove \eqref{Eq:S-N} by induction on $n$. It is trivial to prove \eqref{Eq:S-N} in the case of $n=1$. Now assume that $n>1$ and $k\geqslant n$. Then we have
\begin{align*}
S(N_{k,n})=&\sum\limits_{k_1+\cdots+k_n=k\atop k_i\geqslant 1}(x^{ak_1-1}y+x^{ak_1})S(z_{ak_2}\cdots z_{ak_n})\\
=&\sum\limits_{l=1}^{k-n+1}(z_{al}+x^{al})S(N_{k-l,n-1}).
\end{align*}
By the induction assumption, we obtain
\begin{align*}
S(N_{k,n})=&\sum\limits_{l=1}^{k-n+1}\sum\limits_{i=1}^{n-1}\binom{k-l-i}{k-l-n+1}z_{al}N_{k-l,i}\\
&+\sum\limits_{l=1}^{k-n+1}\sum\limits_{i=1}^{n-1}\binom{k-l-i}{k-l-n+1}x^{al}N_{k-l,i}.
\end{align*}
Since the first term on the right-hand of the above equation is
$$\sum\limits_{i=2}^{n}\sum\limits_{k_1+\cdots+k_{i}=k\atop k_i\geqslant 1}\binom{k-k_1-i+1}{n-i}z_{ak_1}\cdots z_{ak_{i}}$$
and the second term is
$$\sum\limits_{i=1}^{n-1}\sum\limits_{k_1+\cdots+k_{i}=k\atop k_i\geqslant 1}\sum\limits_{l=1}^{k_1-1}\binom{k-l-i}{n-i-1}z_{ak_1}\cdots z_{ak_{i}},$$
we get \eqref{Eq:S-N}.

One can prove \eqref{Eq:S-Inverse-N} similarly, or just from \eqref{Eq:S-N} by using the binomial inversion formula. Finally,  \eqref{Eq:Generating-S-Inverse-N} follows from \eqref{Eq:S-Inverse-N}, and \eqref{Eq:Generating-S-N} follows from \eqref{Eq:S-N}.
\qed

To prove Proposition \ref{Prop:ResSum-E-H}, let us see what equations we need to show in $\mathfrak{h}$.

\begin{lem}
Let $m,n$ be integers with $m\geqslant n\geqslant 1$. Then we have
\begin{align}
\sum\limits_{j=0}^{m-n}(-1)^jz_a^j\ast N_{m-j,n}=(-1)^{m-n}\binom{m}{n}z_a^m
\label{Eq:za-ast-N}
\end{align}
and
\begin{align}
\sum\limits_{j=0}^{m-n}z_a^j\sast N_{m-j,n}=\binom{m}{n}z_a^m.
\label{Eq:za-sast-N}
\end{align}
\end{lem}

\proof We prove \eqref{Eq:za-ast-N}, and one can show \eqref{Eq:za-sast-N} similarly. In fact, \eqref{Eq:za-sast-N} follows from \eqref{Eq:za-ast-N}, because they are equivalent.

We proceed by induction on the pair $(n,m)$ with the lexicographic ordering. If $n=1$, we have to show that for any positive integer $m$, the equation
$$\sum\limits_{j=0}^{m-1}(-1)^jz_a^j\ast z_{a(m-j)}=(-1)^{m-1}mz_a^m$$
holds. While this is just the first equation in (4) of Lemma \ref{Lem:RelSymFunc}. If $m=n$, it is trivial that \eqref{Eq:za-ast-N} holds.

Now assume that $n>1$ and $m>n$. Then we have
\begin{align*}
&\sum\limits_{j=0}^{m-n}(-1)^jz_a^j\ast N_{m-j,n}=N_{m,n}+\sum\limits_{j=1}^{m-n}(-1)^j\sum\limits_{k_1+\cdots+k_n=m-j\atop k_i\geqslant 1}z_a(z_a^{j-1}\ast z_{ak_1}\cdots z_{ak_n})\\
&+\sum\limits_{j=1}^{m-n}(-1)^j\sum\limits_{k_1+\cdots+k_n=m-j\atop k_i\geqslant 1}z_{ak_1}(z_a^{j}\ast z_{ak_2}\cdots z_{ak_n})\\
&+\sum\limits_{j=1}^{m-n}(-1)^j\sum\limits_{k_1+\cdots+k_n=m-j\atop k_i\geqslant 1}z_{a(k_1+1)}(z_a^{j-1}\ast z_{ak_2}\cdots z_{ak_n}).
\end{align*}
Denote the last three terms in the right-hand side of the above equation by $\Sigma_1,\Sigma_2$ and $\Sigma_3$, respectively. By the induction assumption,  we get
$$\Sigma_1=z_a\sum\limits_{j=0}^{m-1-n}(-1)^{j-1}z_a^j\ast N_{m-1-j,n}=(-1)^{m-n}\binom{m-1}{n}z_a^m.$$
And by the induction assumption, we have
\begin{align*}
\Sigma_2=&\sum\limits_{l=1}^{m-n}z_{al}\sum\limits_{j=1}^{m-l-(n-1)}(-1)^jz_a^j\ast N_{m-l-j,n-1}\\
=&\sum\limits_{l=1}^{m-n}z_{al}\left\{(-1)^{m-l-(n-1)}\binom{m-l}{n-1}z_a^{m-l}-N_{m-l,n-1}\right\}\\
=&\sum\limits_{l=1}^{m-n}(-1)^{m-l-(n-1)}\binom{m-l}{n-1}z_{al}z_a^{m-l}-N_{m,n}+z_{a(m-n+1)}z_a^{n-1}.
\end{align*}
Finally, using the induction assumption, we find
\begin{align*}
\Sigma_3=&\sum\limits_{l=2}^{m-n+1}z_{al}\sum\limits_{j=0}^{m-l-(n-1)}(-1)^{j+1}z_a^{j}\ast N_{m-l-j,n-1}\\
=&-\sum\limits_{l=2}^{m-n+1}(-1)^{m-l-(n-1)}\binom{m-l}{n-1}z_{al}z_a^{m-l}.
\end{align*}
Now \eqref{Eq:za-ast-N} follows immediately.
\qed

\noindent {\bf Proof of Proposition \ref{Prop:ResSum-E-H}.}
We get \eqref{Eq:E-ast-F} from \eqref{Eq:za-ast-N}, and get \eqref{Eq:E-sast-F} from \eqref{Eq:za-sast-N}.
Then using Lemma \ref{Lem:RelSymFunc}, we get \eqref{Eq:F-ast-E}.
\qed

Corresponding to \cite[Corollary 1.3]{Hoffman2017}, we have

\begin{cor}
We have
$$S(\mathcal{F}(t,s))\ast \mathcal{F}(t,-s)=1$$
and
$$S^{-1}(\mathcal{F}(t,s))\sast \mathcal{F}(t,-s)=1.$$
\end{cor}

\proof Using Lemma \ref{Lem:RelSymFunc}, Proposition \ref{Prop:ResSum-E-H} and Lemma \ref{Lem:S-F}, we get the results.
\qed

Let $s=1$ and $s=-1$ in Proposition \ref{Prop:ResSum-E-H}, we get the corresponding result of \cite[Corollary 2.2]{Hoffman2017}.

\begin{cor}
We have
$$\mathcal{F}(t,1)=H(t),\quad \mathcal{F}(t,-1)=\overline{H}(-t).$$
Hence for any positive integer $k$, we have
$$\sum\limits_{n=1}^kN_{k,n}=S(z_a^k),\quad \sum\limits_{n=1}^k(-1)^{k+n}N_{k,n}=S^{-1}(z_a^{k}).$$
Moreover, if $a>1$ and $Z_R:\mathfrak{h}^0\rightarrow R$ is a $\mathbb{K}$-linear map, we have
$$\sum\limits_{n=1}^k\sum\limits_{k_1+\cdots+k_n=k\atop k_i\geqslant 1}\zeta_R(ak_1,\ldots,ak_n)=\zeta_R^\star(\{a\}^k)$$
and
$$\sum\limits_{n=1}^k(-1)^{n+k}\sum\limits_{k_1+\cdots+k_n=k\atop k_i\geqslant 1}\zeta_R^{\star}(ak_1,\ldots,ak_n)=\zeta_R(\{a\}^k).$$
\end{cor}

Considering the equations in $\mathfrak{h}^1$ obtained from the two expressions in \eqref{Eq:F-ast-E} of $\mathcal{F}(t,s)$, we get the following result, which is corresponding to \cite[Lemma 2.6]{Hoffman2017}.

\begin{prop}
Let $k,n$ be integers with $k\geqslant n\geqslant 1$. Then we have
\begin{align*}
N_{k,n}=&\sum\limits_{j=0}^{k-n}(-1)^{k-n-j}\binom{k-j}{n}S(z_a^j)\ast z_a^{k-j}\\
=&\sum\limits_{j=0}^{k-n}(-1)^{j}\binom{k-j}{n}S^{-1}(z_a^j)\sast z_a^{k-j}.
\end{align*}
Moreover, if $a>1$ and $Z_R:\mathfrak{h}^0_{\ast}\rightarrow R$ is an algebra homomorphism, we have
\begin{align*}
&\sum\limits_{k_1+\cdots+k_n=k\atop k_i\geqslant 1}\zeta_R(ak_1,\ldots,ak_n)=\sum\limits_{j=0}^{k-n}(-1)^{k-n-j}\binom{k-j}{n}\zeta_R^{\star}(\{a\}^j)\zeta_R(\{a\}^{k-j}),\\
&\sum\limits_{k_1+\cdots+k_n=k\atop k_i\geqslant 1}\zeta_R^{\star}(ak_1,\ldots,ak_n)=\sum\limits_{j=0}^{k-n}(-1)^{j}\binom{k-j}{n}\zeta_R(\{a\}^j)\zeta_R^{\star}(\{a\}^{k-j}).
\end{align*}
\end{prop}

If $a$ is even, using Theorem \ref{Thm:Ev-2k}, we get the following restricted sum formulas.

\begin{thm}\label{Thm:ResSum-2a}
Let $m$ be a positive integer and the $\mathbb{K}$-linear map $Z_R:\mathfrak{h}^0\rightarrow R$ satisfy the regularized double shuffle relations. Then for any integers $k,n$ with $k\geqslant n\geqslant 1$, we have
\begin{align}
&\sum\limits_{k_1+\cdots+k_n=k\atop k_i\geqslant 1}\zeta_R(2mk_1,\ldots,2mk_n)=(-1)^n\sum\limits_{j=0}^{k-n}\binom{k-j}{n}\sum\limits_{{n_0+\cdots+n_{m-1}=mj\atop l_0+\cdots+l_{m-1}=m(k-j)}\atop n_i,l_i\geqslant 0}\nonumber\\
&\qquad \times \left[\prod\limits_{i=0}^{m-1}\frac{(2-4^{n_i})B_{2n_i}}{(2n_i)!(2l_i+1)!}\right]
\rho_m^{\sum\limits_{i=0}^{m-1}2i(n_i+l_i)}\frac{\lambda^{2km}}{4^{km}}
\label{Eq:Res-Hoffman-MZV-general}
\end{align}
and
\begin{align}
&\sum\limits_{k_1+\cdots+k_n=k\atop k_i\geqslant 1}\zeta_R^{\star}(2mk_1,\ldots,2mk_n)=\sum\limits_{j=n}^{k}\binom{j}{n}\sum\limits_{{n_0+\cdots+n_{m-1}=mj\atop l_0+\cdots+l_{m-1}=m(k-j)}\atop n_i,l_i\geqslant 0}\nonumber\\
&\qquad \times \left[\prod\limits_{i=0}^{m-1}\frac{(2-4^{n_i})B_{2n_i}}{(2n_i)!(2l_i+1)!}\right]
\rho_m^{\sum\limits_{i=0}^{m-1}2i(n_i+l_i)}\frac{\lambda^{2km}}{4^{km}},
\label{Eq:Res-Hoffman-MZSV-general}
\end{align}
where $\rho_m=e^{\frac{\pi \sqrt{-1}}{m}}$.
\end{thm}

When $Z_R=Z$ and $\lambda=2\pi\sqrt{-1}$, \eqref{Eq:Res-Hoffman-MZV-general} is much simpler than \cite[Theorem 7.1]{Komori-Matsumoto-Tsumura}.

Let $m=1$ in Theorem \ref{Thm:ResSum-2a}, we have
\begin{align}
\sum\limits_{k_1+\cdots+k_n=k\atop k_i\geqslant 1}\zeta_R(2k_1,\ldots,2k_n)=&\frac{(-1)^n\lambda^{2k}}{4^k(2k+1)!}\sum\limits_{j=0}^{k-n}\binom{k-j}{n}\binom{2k+1}{2j}(2-4^{j})B_{2j}
\label{Eq:Res-Hoffman}
\end{align}
and
\begin{align*}
\sum\limits_{k_1+\cdots+k_n=k\atop k_i\geqslant 1}\zeta_R^{\star}(2k_1,\ldots,2k_n)=&\frac{\lambda^{2k}}{4^k(2k+1)!}\sum\limits_{j=n}^{k}\binom{j}{n}\binom{2k+1}{2j}(2-4^{j})B_{2j},
\end{align*}
which correspond to \cite[Theorem 1.4]{Hoffman2017} and \cite[Eq. (23)]{Gencev} in the case of $Z_R=Z$ and $\lambda=2\pi\sqrt{-1}$, respectively.
Let $m=2$ in Theorem \ref{Thm:ResSum-2a}, and use \eqref{Eq:z-4-n-new} instead of \eqref{Eq:z-2k-n}, we get
\begin{align*}
\sum\limits_{k_1+\cdots+k_n=k\atop k_i\geqslant 1}\zeta_R(4k_1,\ldots,4k_n)=&(-1)^{k+n}\frac{2\lambda^{4k}}{4^k(4k+2)!}\sum\limits_{j=0}^{k-n}\frac{(-1)^j}{4^j}\binom{k-j}{n}\binom{4k+2}{4j}\\
&\times \sum\limits_{i=0}^{2j}(-1)^i\binom{4j}{2i}(2-4^i)(2-4^{2j-i})B_{2i}B_{4j-2i}
\end{align*}
and
\begin{align*}
\sum\limits_{k_1+\cdots+k_n=k\atop k_i\geqslant 1}\zeta_R^{\star}(4k_1,\ldots,4k_n)=&(-1)^{k}\frac{2\lambda^{4k}}{4^k(4k+2)!}\sum\limits_{j=n}^{k}\frac{(-1)^j}{4^j}\binom{j}{n}\binom{4k+2}{4j}\\
&\times \sum\limits_{i=0}^{2j}(-1)^i\binom{4j}{2i}(2-4^i)(2-4^{2j-i})B_{2i}B_{4j-2i}.
\end{align*}
When $Z_R=Z$ and $\lambda=2\pi\sqrt{-1}$, the first equation above is just \cite[Eq. (20)]{Gencev}. Hence we give a pure algebraic proof of \cite[Conjecture 4.1]{Gencev} (compare with \cite[Theorem 1.1]{Yuan-Zhao}).

\begin{rem}
In \cite{Chen-Chung-Eie}, K.-W. Chen, C.-L. Chung and M. Eie also considered the restricted sum
$$\sum\limits_{k_1+\cdots+k_n=k\atop k_i\geqslant 1}\zeta(ak_1,\ldots,ak_n).$$
\end{rem}

We would like to state below the corresponding results of \cite[Proposition 2.4 and Proposition 2.5]{Hoffman2017}, although it will not be used.

\begin{prop}
\begin{itemize}
  \item[(1)] We have
  $$t\frac{\partial\mathcal{F}}{\partial t}(t,s)+(1-s)\frac{\partial\mathcal{F}}{\partial s}(t,s)=tP(t)\ast \mathcal{F}(t,s)$$
and
$$t\frac{\partial\mathcal{F}}{\partial t}(t,s)-(s+1)\frac{\partial\mathcal{F}}{\partial s}(t,s)=-tP(t)\sast \mathcal{F}(t,s);$$
  \item[(2)] For any positive integers $k,n$ with $k\geqslant n+1$, we have
$$\sum\limits_{j=1}^{k-n}z_{aj}\ast N_{k-j,n}=(k-n)N_{k,n}+(n+1)N_{k,n+1}$$
and
$$\sum\limits_{j=1}^{k-n}z_{aj}\sast N_{k-j,n}=(n-k)N_{k,n}+(n+1)N_{k,n+1}.$$
\end{itemize}
\end{prop}

\proof
We get (1) from Proposition \ref{Prop:ResSum-E-H} and Lemma \ref{Lem:RelSymFunc}. Then (2) is deduced from (1). In fact, (1) and (2) are equivalent.
\qed

Below we set $a=2$. Then we get the result corresponding to \cite[Theorem 1.2]{Hoffman2017}.

\begin{prop}
If $a=2$ and the $\mathbb{K}$-linear map $Z_R:\mathfrak{h}^0\rightarrow R$ satisfies the regularized double shuffle relations, we have
\begin{align*}
&Z_R(\mathcal{F}(t,s))=\frac{e^{\lambda\sqrt{(1-s)t}/2}-e^{-\lambda\sqrt{(1-s)t}/2}}{\sqrt{1-s}\left(e^{\lambda\sqrt{t}/2}-e^{-\lambda\sqrt{t}/2}\right)},\\
&Z_R^{\star}(\mathcal{F}(t,s))=\frac{\sqrt{1+s}\left(e^{\lambda\sqrt{t}/2}-e^{-\lambda\sqrt{t}/2}\right)}{e^{\lambda\sqrt{(1+s)t}/2}-e^{-\lambda\sqrt{(1+s)t}/2}}.
\end{align*}
\end{prop}

\proof If $a=2$, by \eqref{Eq:MZVk-MZSVk} and \eqref{Eq:Generating-2-n}, we have
\begin{align}
&Z_R(E((s-1)t))=\frac{e^{\lambda\sqrt{(1-s)t}/2}-e^{-\lambda\sqrt{(1-s)t}/2}}{\lambda\sqrt{(1-s)t}},\label{Eq:Z-E}\\
&Z_R(H(t))=\frac{\lambda\sqrt{t}}{e^{\lambda\sqrt{t}/2}-e^{-\lambda\sqrt{t}/2}},\label{Eq:Z-H}\\
&Z_R^{\star}(E((s+1)t))=\frac{\lambda\sqrt{(1+s)t}}{e^{\lambda\sqrt{(1+s)t}/2}-e^{-\lambda\sqrt{(1+s)t}/2}},\nonumber\\
&Z_R^{\star}(\overline{H}(-t))=\frac{e^{\lambda\sqrt{t}/2}-e^{-\lambda\sqrt{t}/2}}{\lambda\sqrt{t}}.\nonumber
\end{align}
Then we get the result from Proposition \ref{Prop:ResSum-E-H}.
\qed

Finally, corresponding to \cite[Theorem 1.1]{Hoffman2017}, we have

\begin{thm}\label{Thm:res-Hoffman-two}
Let $k,n$ be integers with $k\geqslant n\geqslant 1$ and the $\mathbb{K}$-linear map $Z_R:\mathfrak{h}^0\rightarrow R$ satisfy the regularized double shuffle relations. Then we have
\begin{align*}
&\sum\limits_{k_1+\cdots+k_n=k\atop k_i\geqslant 1}\zeta_R(2k_1,\ldots,2k_n)=\sum\limits_{j=0}^{\left[\frac{n-1}{2}\right]}\frac{\lambda^{2j}}{2^{2n-2}(2j+1)!}\binom{2n-2j-1}{n}\zeta_R(2k-2j)\\
=&\frac{1}{2^{2n-2}}\binom{2n-1}{n}\zeta_R(2k)-\sum\limits_{j=1}^{\left[\frac{n-1}{2}\right]}\frac{\binom{2n-2j-1}{n}}{2^{2n-3}(2j+1)B_{2j}}\zeta_R(2j)\zeta_R(2k-2j).
\end{align*}
\end{thm}

One can prove Theorem \ref{Thm:res-Hoffman-two} by using \eqref{Eq:Res-Hoffman} and \cite[Theorem 1.5]{Hoffman2017}. Here we provide a direct proof, which is similar as the proof of \cite[Theorem 1.5]{Hoffman2017}. In fact, by \eqref{Eq:Z-E}, we have
\begin{align*}
Z_R(E((s-1)t))=&\sum\limits_{n=0}^\infty\frac{\lambda^{2n}(1-s)^nt^n}{2^{2n}(2n+1)!}\\
=&\sum\limits_{m=0}^\infty(-1)^m\left(\sum\limits_{n\geqslant m}\binom{n}{m}\frac{\lambda^{2n}}{2^{2n}(2n+1)!}t^n\right)s^m.
\end{align*}
Set
$$f(t)=\frac{e^{\lambda\sqrt{t}/2}-e^{-\lambda\sqrt{t}/2}}{\lambda\sqrt{t}}=\sum\limits_{n=0}^\infty\frac{\lambda^{2n}}{2^{2n}(2n+1)!}t^n,$$
Then by \eqref{Eq:Z-H}, we have
$$Z_R(\mathcal{F}(t,s))=\sum\limits_{m=0}^\infty G_m(t)s^m,$$
where
$$G_m(t)=(-1)^m\frac{1}{f(t)}\frac{t^m}{m!}\frac{d^m}{dt^m}f(t).$$
Set
$$g(t)=e^{\lambda\sqrt{t}/2}+e^{-\lambda\sqrt{t}/2}.$$

\begin{lem}\label{Lem:D-f}
For any nonnegative integer $m$, we have
\begin{align*}
\frac{(-1)^m2^{2m}t^m}{m!}\frac{d^m}{dt^m}f(t)=&\left(\sum\limits_{j=0}^{\left[\frac{m}{2}\right]}\frac{\lambda^{2j}}{(2j)!}\binom{2m-2j}{m}t^j\right)f(t)\\
&-\left(\sum\limits_{j=0}^{\left[\frac{m-1}{2}\right]}\frac{\lambda^{2j}}{(2j+1)!}\binom{2m-2j-1}{m}t^j\right)g(t).
\end{align*}
\end{lem}

\proof
We proceed by induction on $m$. It is trivial to check the case of $m=0$. Now assume that the formula is valid for $m$, then by
$$\frac{d}{dt}f(t)=-\frac{1}{2t}f(t)+\frac{1}{4t}g(t),\quad \frac{d}{dt}g(t)=\frac{\lambda^2}{4}f(t),$$
we get
$$\frac{(-1)^{m+1}2^{2(m+1)}t^{m+1}}{(m+1)!}\frac{d^{m+1}}{dt^{m+1}}f(t)=A_mf(t)+B_mg(t),$$
where
\begin{align*}
A_m=&-\frac{4}{m+1}\sum\limits_{j=0}^{\left[\frac{m}{2}\right]}\frac{\lambda^{2j}}{(2j)!}\binom{2m-2j}{m}(j-m)t^j+\frac{2}{m+1}\sum\limits_{j=0}^{\left[\frac{m}{2}\right]}\frac{\lambda^{2j}}{(2j)!}\binom{2m-2j}{m}t^j\\
&+\frac{t}{m+1}\sum\limits_{j=0}^{\left[\frac{m-1}{2}\right]}\frac{\lambda^{2j+2}}{(2j+1)!}\binom{2m-2j-1}{m}t^j,\\
B_m=&\frac{4}{m+1}\sum\limits_{j=0}^{\left[\frac{m-1}{2}\right]}\frac{\lambda^{2j}}{(2j+1)!}\binom{2m-2j-1}{m}(j-m)t^j\\
&-\frac{1}{m+1}\sum\limits_{j=0}^{\left[\frac{m}{2}\right]}\frac{\lambda^{2j}}{(2j)!}\binom{2m-2j}{m}t^j.
\end{align*}
We  have
\begin{align*}
A_m=&\frac{2}{m+1}\sum\limits_{j=0}^{\left[\frac{m}{2}\right]}\frac{\lambda^{2j}}{(2j)!}\binom{2m-2j}{m}(2m-2j+1)t^j\\
&+\frac{1}{m+1}\sum\limits_{j=1}^{\left[\frac{m+1}{2}\right]}\frac{\lambda^{2j}}{(2j-1)!}\binom{2m-2j+1}{m}t^j\\
=&\sum\limits_{j=0}^{\left[\frac{m}{2}\right]}\frac{\lambda^{2j}}{(2j)!}\binom{2m-2j+1}{m+1}t^j\\
&+\frac{1}{m+1}\sum\limits_{j=0}^{\left[\frac{m}{2}\right]}\frac{\lambda^{2j}}{(2j)!}\binom{2m-2j+1}{m}(m-2j+1)t^j\\
&+\frac{1}{m+1}\sum\limits_{j=1}^{\left[\frac{m+1}{2}\right]}\frac{\lambda^{2j}}{(2j-1)!}\binom{2m-2j+1}{m}t^j\\
=&\sum\limits_{j=0}^{\left[\frac{m}{2}\right]}\frac{\lambda^{2j}}{(2j)!}\binom{2m-2j+1}{m+1}t^j
+\sum\limits_{j=0}^{\left[\frac{m}{2}\right]}\frac{\lambda^{2j}}{(2j)!}\binom{2m-2j+1}{m}t^j\\
&-\frac{1}{m+1}\sum\limits_{j=1}^{\left[\frac{m}{2}\right]}\frac{\lambda^{2j}}{(2j-1)!}\binom{2m-2j+1}{m}t^j\\
&+\frac{1}{m+1}\sum\limits_{j=1}^{\left[\frac{m+1}{2}\right]}\frac{\lambda^{2j}}{(2j-1)!}\binom{2m-2j+1}{m}t^j\\
=&\sum\limits_{j=0}^{\left[\frac{m}{2}\right]}\frac{\lambda^{2j}}{(2j)!}\binom{2m-2j+2}{m+1}t^j+C_m,
\end{align*}
where
$$C_m=\begin{cases}
0 & \text{if\;} m\text{\;is even},\\
\frac{\lambda^{m+1}}{(m+1)!}t^{\frac{m+1}{2}} & \text{if\;} m\text{\;is odd},
\end{cases}$$
which implies that
$$A_m=\sum\limits_{j=0}^{\left[\frac{m+1}{2}\right]}\frac{\lambda^{2j}}{(2j)!}\binom{2m-2j+2}{m+1}t^j.$$
Similarly, we have
$$B_m=-\sum\limits_{j=0}^{\left[\frac{m}{2}\right]}\frac{\lambda^{2j}}{(2j+1)!}\binom{2m-2j+1}{m+1}t^j.$$
Hence we finish the proof.
\qed

By Lemma \ref{Lem:D-f}, we get
\begin{align*}
G_m(t)=&\frac{1}{2^{2m}}\sum\limits_{j=0}^{\left[\frac{m}{2}\right]}\frac{\lambda^{2j}}{(2j)!}\binom{2m-2j}{m}t^j\\
&-\frac{1}{2^{2m}}\left(\sum\limits_{j=0}^{\left[\frac{m-1}{2}\right]}\frac{\lambda^{2j}}{(2j+1)!}\binom{2m-2j-1}{m}t^j\right)\frac{g(t)}{f(t)}.
\end{align*}

\begin{lem}\label{Lem:g-div-f}
Let the $\mathbb{K}$-linear map $Z_R:\mathfrak{h}^0\rightarrow R$ satisfy the regularized double shuffle relations. Then we have
$$\frac{g(t)}{f(t)}=-4\sum\limits_{n=0}^\infty\zeta_R(2n)t^n.$$
\end{lem}

\proof
Applying the operator $\frac{d}{dt}$ to both sides of the equation
$$e^{\lambda\sqrt{t}/2}-e^{-\lambda\sqrt{t}/2}=\lambda\sqrt{t}\exp\left(-\sum\limits_{n=1}^\infty\frac{1}{n}\zeta_R(2n)t^n\right),$$
we get
$$\frac{\lambda}{4\sqrt{t}}\left(e^{\lambda\sqrt{t}/2}+e^{-\lambda\sqrt{t}/2}\right)=\left(e^{\lambda\sqrt{t}/2}-e^{-\lambda\sqrt{t}/2}\right)
\left(\frac{1}{2t}-\sum\limits_{n=1}^\infty\zeta_R(2n)t^{n-1}\right).$$
Then the result follows from $\zeta_R(0)=-\frac{1}{2}$.
\qed

Using Lemma \ref{Lem:g-div-f}, we have
\begin{align*}
G_m(t)=&\frac{1}{2^{2m}}\sum\limits_{j=0}^{\left[\frac{m}{2}\right]}\frac{\lambda^{2j}}{(2j)!}\binom{2m-2j}{m}t^j\\
&+\frac{1}{2^{2m-2}}\sum\limits_{0\leqslant j\leqslant\frac{m-1}{2}\atop k\geqslant j}\frac{\lambda^{2j}}{(2j+1)!}\binom{2m-2j-1}{m}\zeta_R(2k-2j)t^k.
\end{align*}
Finally we get
\begin{align*}
&Z_R(\mathcal{F}(t,s))=\sum\limits_{m\geqslant 0\atop 0\leqslant j\leqslant \frac{m}{2}}\frac{1}{2^{2m}}\frac{\lambda^{2j}}{(2j)!}\binom{2m-2j}{m}t^js^m\\
&\quad+\sum\limits_{m\geqslant 0,0\leqslant j\leqslant\frac{m-1}{2}\atop k\geqslant j}\frac{1}{2^{2m-2}}\frac{\lambda^{2j}}{(2j+1)!}\binom{2m-2j-1}{m}\zeta_R(2k-2j)t^ks^m.
\end{align*}
Comparing the coefficients of $t^ks^n$, we get Theorem \ref{Thm:res-Hoffman-two}.\qed


\subsection{Some evidence of Brown-Zagier relation}

In \cite{Zagier2012}, D. Zagier proved that for any nonnegative integers $n$ and $m$, we have
\begin{align}
&\zeta(\{2\}^n,3,\{2\}^m)=\sum\limits_{r=1}^{m+n+1}
c^r_{m,n}\zeta(2r+1)\zeta(\{2\}^{m+n+1-r}),
\label{Eq:Z-2-3-2}\\
&\zeta^{\star}(\{2\}^n,3,\{2\}^m)=\sum\limits_{r=1}^{m+n+1}c_{m,n}^{\star,r}\zeta(2r+1)
\zeta^{\star}(\{2\}^{m+n+1-r}),
\label{Eq:SZ-2-3-2}
\end{align}
where
\begin{align*}
&c^r_{m,n}=2(-1)^r\left\{\binom{2r}{2m+2}-(1-2^{-2r})\binom{2r}{2n+1}\right\},\\
&c^{\star,r}_{m,n}=-2\left\{\binom{2r}{2m}-\delta_{r,m}-(1-2^{-2r})\binom{2r}{2n+1}\right\}.
\end{align*}
The formula \eqref{Eq:Z-2-3-2} was used in \cite{Brown}. The first named author gave a simpler proof of this evaluation formula in \cite{Li2013-3}, and  T. Terasoma showed that this formula is a consequence of the associator relations in \cite{Terasoma-2}. We give the following conjecture.

\begin{conj}
Let the $\mathbb{K}$-linear map $Z_R:\mathfrak{h}^0\rightarrow R$ satisfy the regularized double shuffle relations. Then for any nonnegative integers $m,n$, we have
\begin{align}
&\zeta_R(\{2\}^n,3,\{2\}^m)=\sum\limits_{r=1}^{m+n+1}
c^r_{m,n}\zeta_R(2r+1)\zeta_R(\{2\}^{m+n+1-r}),
\label{Eq:Brown-Zagier-MZV}\\
&\zeta^{\star}_R(\{2\}^n,3,\{2\}^m)=\sum\limits_{r=1}^{m+n+1}c_{m,n}^{\star,r}\zeta_R(2r+1)
\zeta^{\star}_R(\{2\}^{m+n+1-r}).
\label{Eq:Brown-Zagier-MZSV}
\end{align}
\end{conj}

In \cite{Li2013-2}, the first named author showed that \eqref{Eq:Brown-Zagier-MZV} and \eqref{Eq:Brown-Zagier-MZSV} are equivalent. Here we show that under the stuffle products, summing over $n+m=k$ for some fixed $k$ of both sides of \eqref{Eq:Brown-Zagier-MZV} (resp. \eqref{Eq:Brown-Zagier-MZSV}), we indeed get equality.

\begin{prop}
For any positive integers $a,b$ and any nonnegative integer $k$, we have
\begin{align*}
\sum\limits_{n+m=k\atop n,m\geqslant 0}\sum\limits_{r=1}^{m+n+1}c_{m,n}^rz_{b+(r-1)a}\ast z_a^{m+n+1-r}=&\sum\limits_{m=0}^{k}(-1)^mz_{b+ma}\ast z_{a}^{k-m}\\
=&\sum\limits_{n=0}^kz_a^nz_bz_a^{k-n},
\end{align*}
and
\begin{align*}
\sum\limits_{n+m=k\atop n,m\geqslant 0}\sum\limits_{r=1}^{m+n+1}c_{m,n}^{\star,r}z_{b+(r-1)a}\sast z_a^{m+n+1-r}=\sum\limits_{m=0}^{k}z_{b+ma}\sast z_{a}^{k-m}=\sum\limits_{n=0}^kz_a^nz_bz_a^{k-n}.
\end{align*}
Moreover, if $Z_R:\mathfrak{h}^0_{\ast}\rightarrow R$ is an algebra homomorphism, we have
$$\sum\limits_{n+m=k\atop n,m\geqslant 0}\sum\limits_{r=1}^{m+n+1}
c^r_{m,n}\zeta_{R}(2r+1)\zeta_{R}(\{2\}^{m+n+1-r})=\sum\limits_{n+m=k\atop n,m\geqslant 0}\zeta_{R}(\{2\}^n,3,\{2\}^m)$$
and
$$\sum\limits_{n+m=k\atop n,m\geqslant 0}\sum\limits_{r=1}^{m+n+1}
c^{\star,r}_{m,n}\zeta_{R}(2r+1)\zeta_{R}^{\star}(\{2\}^{m+n+1-r})=\sum\limits_{n+m=k\atop n,m\geqslant 0}\zeta_{R}^{\star}(\{2\}^n,3,\{2\}^m).$$
\end{prop}

\proof
We have
\begin{align*}
&\sum\limits_{n+m=k\atop n,m\geqslant 0}\sum\limits_{r=1}^{m+n+1}c_{m,n}^rz_{b+(r-1)a}\ast z_a^{m+n+1-r}
=\sum\limits_{n+m=k\atop n,m\geqslant 0}\sum\limits_{r=1}^{k+1}c_{m,n}^rz_{b+(r-1)a}\ast z_a^{k+1-r}\\
=&\sum\limits_{r=1}^{k+1}\left(\sum\limits_{n+m=k\atop n,m\geqslant 0}c_{m,n}^r\right)z_{b+(r-1)a}\ast z_a^{k+1-r}
\end{align*}
and
\begin{align*}
&\sum\limits_{n+m=k\atop n,m\geqslant 0}\sum\limits_{r=1}^{m+n+1}c_{m,n}^{\star,r}z_{b+(r-1)a}\sast z_a^{m+n+1-r}=\sum\limits_{r=1}^{k+1}\left(\sum\limits_{n+m=k\atop n,m\geqslant 0}c_{m,n}^{\star,r}\right)z_{b+(r-1)a}\sast z_a^{k+1-r}.
\end{align*}

Now for any integer $r$ with $1\leqslant r\leqslant k+1$, we have
\begin{align*}
\sum\limits_{n+m=k\atop n,m\geqslant 0}c^r_{m,n}=&2(-1)^r\sum\limits_{n+m=k\atop n,m\geqslant 0}\binom{2r}{2m+2}-2(-1)^r(1-2^{-2r})
\sum\limits_{n+m=k\atop m,n\geqslant 0}\binom{2r}{2n+1}\\
=&2(-1)^r\sum\limits_{i=1}^r\binom{2r}{2i}-2(-1)^r(1-2^{-2r})\sum\limits_{i=1}^r\binom{2r}{2i-1}
\end{align*}
and
\begin{align*}
\sum\limits_{n+m=k\atop n,m\geqslant 0}c^{\star,r}_{m,n}=&-2\sum\limits_{i=0}^r\binom{2r}{2i}+2+2(1-2^{-2r})\sum\limits_{i=1}^r\binom{2r}{2i-1}.
\end{align*}
Since
$$2\sum\limits_{i=0}^r\binom{2r}{2i}=2^{2r}=2\sum\limits_{i=1}^r\binom{2r}{2i-1},$$
we get
$$
\sum\limits_{n+m=k\atop n,m\geqslant 0}c^r_{m,n}=(-1)^{r+1},\qquad
\sum\limits_{n+m=k\atop n,m\geqslant 0}c^{\star,r}_{m,n}=1.$$
Then with the help of Lemma \ref{Lem:a-b-ast} we get the result.
\qed


\section{Further remarks}\label{Sec:Remarks}

Besides the regularized double shuffle relations, there are some other systems of relations of multiple zeta values, which are also conjectured to give all algebraic relations. For example, the following relations are known
\begin{itemize}
  \item associator relations (\cite{Drinfeld,Furusho2010});
  \item Kawashima's relation (\cite{Kawashima,Tanaka});
  \item Kaneko-Yamamoto relation (\cite{Kaneko-Yamamoto}).
\end{itemize}
We recall some relationship between the regularized double shuffle relations and the relations listed above. It was proved in \cite{Deligne-Terasoma,Furusho2011} that the regularized double shuffle relations can be deduced from the associator relations.  In \cite{Kaneko-Yamamoto}, it was proved that the Kaneko-Yamamoto relations are equivalent to the regularized double shuffle relations, and under the duality formula, the Kawashima relation can be deduced from the Kaneko-Yamamoto relations. Since all the four systems of relations should be equivalent, it is interesting to consider other implication relations among these relations.


\end{document}